\documentclass[a4paper,UKenglish]{lipics_md}
\usepackage{microtype}

\title{Expected Sizes of Poisson--Delaunay Mosaics and Their
       Discrete Morse Functions\footnote{This work is partially
         supported by the {\sc Toposys} project FP7-ICT-318493-STREP,
         by ESF under the ACAT Research Network Programme,
         and by the FWF and DFG within the SFB-Transregio Programme 109
         in Discrete Differential Geometry, grant no.\ I2979.}}
\titlerunning{Discrete Delaunay--Morse Function}

\author[1]{Herbert Edelsbrunner}
\author[1]{Anton Nikitenko}
\author[2]{Matthias Reitzner}
\affil[1]{IST Austria (Institute of Science and Technology Austria),
  Am Campus 1, \\ 3400 Klosterneuburg, Austria,
  \texttt{edels@ist.ac.at}, \texttt{anton.nikitenko@ist.ac.at}}
\affil[2]{Mathematics Department, University of Osnabr\"{u}ck,
  Albrechtstra{\ss}e 28a, \\ 49076 Osnabr\"{u}ck, Germany, 
  \texttt{mreitzne@uni-osnabrueck.de}}
\authorrunning{H. Edelsbrunner, A. Nikitenko and M. Reitzner}

\subjclass{I.3.5 Computational Geometry and Object Modeling, G.3 Probability and Statistics, G.2 Discrete Mathematics}
\mathsubjclass{60D05 Geometric probability and stochastic geometry, 68U05 Computer graphics; computational geometry} 
\keywords{Poisson point process, Delaunay mosaic, discrete Morse theory,
   critical simplices, intervals, stochastic geometry, integral geometry, 
typical simplex.}

\newcommand{\mm}[1] {\ifmmode{#1}\else{\mbox{\(#1\)}}\fi}

\newcommand{\denselist}{\itemsep 0pt\parsep=1pt\partopsep 0pt}

\newcommand{\ignore}[1]{}

\newcommand{\ourproof}{\begin{proof}}
\newcommand{\eop}{\end{proof}}  %

\newcommand{\Aspace}        {\mm{{\mathbb A}}}
\newcommand{\Bspace}        {\mm{{\mathbb B}}}
\newcommand{\Hspace}        {\mm{{\mathbb H}}}
\newcommand{\RPspace}       {\mm{{\mathbb R}{\mathbb P}}}
\newcommand{\Rspace}        {\mm{{\mathbb R}}}
\newcommand{\Sspace}        {\mm{{\mathbb S}}}

\newcommand{\Zspace}        {\mm{{\mathbb Z}}}
\newcommand{\EE}            {\mm{{\mathbb E}}}
\newcommand{\Expected}[1]   {\mm{{\EE}{[{#1}]}}}
\newcommand{\PP}            {\mm{{\mathbb P}}}
\newcommand{\Probable}[2]   {\mm{{\PP}_{#1}{[{#2}]}}}
\newcommand{\Pempty}[1]     {\mm{{\PP}_{\emptyset}{[{#1}]}}}
\newcommand{\One}[1]        {\mm{{{\bf 1}_{#1}}}}

\newcommand{\Region}        {\mm{\Omega}}

\newcommand{\ttt}           {\mm{{\bf t}}}
\newcommand{\Minus}[1]      {\mm{{\#}{#1}}}
\newcommand{\uuu}           {\mm{{\bf u}}}

\newcommand{\xxx}           {\mm{{\bf x}}}
\newcommand{\LGrass}[2]     {\mm{{\cal L}_{#1}^{#2}}}

\newcommand{\Delaunay}[1]   {\mm{{\rm Del\,}{{#1}}}}

\newcommand{\Rf}            {\mm{{\rm Rad}}}
\newcommand{\density}       {\mm{\rho}}
\newcommand{\littleoh}[1]   {\mm{o{\left({#1}\right)}}}
\newcommand{\bigoh}[1]      {\mm{O{({#1})}}}

\newcommand{\Ccon}[3]       {\mm{{C}_{{#1},{#2}}^{#3}}}
\newcommand{\Dcon}[2]       {\mm{{D}_{#1}^{#2}}}
\newcommand{\ccon}[3]       {\mm{{c}_{{#1},{#2}}^{#3}}}
\newcommand{\dcon}[2]       {\mm{{d}_{#1}^{#2}}}

\newcommand{\Factor}[2]     {\mm{{\rm Factor}}{({#1},{#2})}}
\newcommand{\FactorOnly}    {\mm{{\rm Factor}}}
\newcommand{\xMnt}[3]       {\mm{{\rm Mnt}_{2}{({{#1};{#2},{#3}})}}}
\newcommand{\xMntOnly}      {\mm{{\rm Mnt}_{2}}}
\newcommand{\Mnt}[3]        {\mm{{\rm Mnt}_{1}{({{#1},{#2};{#3}})}}}
\newcommand{\MntOnly}       {\mm{{\rm Mnt}_{1}}}
\newcommand{\EAux}[3]       {\mm{{E}_{{#1},{#2}}^{#3}}}
\newcommand{\Vsign}[1]      {\mm{{V}_{#1}}}
\newcommand{\Vcone}[2]      {\mm{{V}_{#1}^{#2}}}
\newcommand{\Volume}[1]     {\mm{\rm Vol}{({#1})}}
\newcommand{\Voronoi}[1]    {\mm{{\rm Vor}{({#1})}}}
\newcommand{\Signum}[1]     {\mm{{\rm sgn}{({#1})}}}

\newcommand{\card}[1]       {\mm{|{#1}|}}
\newcommand{\dime}[1]       {\mm{\rm dim\,}{#1}}

\newcommand{\diff}          {\mm{\rm \,d}}
\newcommand{\norm}[1]       {\mm{\|{#1}\|}}
\newcommand{\Edist}[2]      {\mm{\|{#1}-{#2}\|}}
\newcommand{\dd}            {\mm{\delta}}
\newcommand{\ee}            {\mm{\varepsilon}}

\newcommand{\ourparagraph}[1]  {\vspace{0.1in} \noindent \textbf{#1}}

\begin{document}
\maketitle

\begin{abstract}
  Mapping every simplex in the Delaunay mosaic of a discrete point set
  to the radius of the smallest empty circumsphere gives a generalized
  discrete Morse function.
  Choosing the points from a Poisson point process in $\Rspace^n$,
  we study the expected number of simplices in the
  Delaunay mosaic as well as the expected number of
  critical simplices and non-singular intervals
  in the corresponding generalized discrete gradient.
  Observing connections with other probabilistic models,
  we obtain precise expressions for the expected numbers in low dimensions.
  In particular, we get the expected numbers of simplices in the
  Poisson--Delaunay mosaic in dimensions $n \leq 4$.
\end{abstract}

\section{Introduction}
\label{sec:1}

One motivation for the work reported in this paper is the desire
to reconstruct surfaces from point sets;
see \cite{Dey11} and in particular the Wrap algorithm described in \cite{Ede03}.
While the points usually describe a distinctive shape and are
therefore not random, they are affected by noise and display
random features locally.
To effectively cope with local noise is a necessary component of every
high quality surface reconstruction software.
Another motivation derives from the work in topological data analysis;
see \cite{Car09}.
The understanding of random data sets provides the necessary 
background against which we can assess possibly non-random features.
The most natural and accessible model in this direction is to choose
points randomly according to a Poisson point process. 
This set of random points gives rise to the Poisson--Delaunay mosaic,
which generically is a simplicial complex tessellating the space.
This tessellation is a well investigated classical mosaic;
see e.g.\ Miles \cite{Mil70,Mil71} and Chapter 10 in the book
of Schneider and Weil \cite{ScWe08} on stochastic geometry.
Yet many questions are still open.
For example, the expected number of simplices in a set of volume one
was only known in dimensions two and three.

More interesting from a topological point of view is the \emph{Alpha complex}
or \emph{Delaunay complex} for radius $r \geq 0$ --- the subcomplex
consisting of all Delaunay simplices with circumradius at most $r$;
see \cite{EdMu94}.
Its simplices do not cover the entire Euclidean space and can therefore
form cycles are other topological features.
Random simplicial complexes have recently found considerable attention.
In particular, topological characteristics of \v{C}ech and Rips complexes
over Poisson point processes have been investigated
in work of Kahle \cite{Kah11,Kah14}, Bobrowski and Weinberger \cite{BoWe15},
Bobrowski and Adler \cite{BoAd14}, and Decreusefond et al.\ \cite{DFRV14}.
We add to this direction by calculating the expected number of simplices in
Delaunay complexes as well as the expected number of critical simplices
and non-singular intervals of the corresponding radius function.
In an equivalent formulation, this yields the distribution function of the 
circumradius of the typical $j$-dimensional simplex of the Poisson--Delaunay mosaic.
The distribution of the circumradius for the top-dimensional case is due 
to Miles \cite{Mil74}, see also M\o ller \cite{Mol89},
but the general cases have been unknown so far. 

\ourparagraph{Results.}
We introduce a few concepts to formally state our results.
Letting $X$ be sampled from a Poisson point process with density
$\density > 0$ in $\Rspace^n$, we consider the \emph{Delaunay triangulation}
or \emph{Delaunay mosaic}, $\Delaunay{X}$, which consists of simplices of 
dimensions $0, 1, \ldots, n$. 
To each simplex we assign the smallest ball containing no points of $X$
in the interior and all vertices of the simplex on its boundary.
We denote by  $\Rf \colon \Delaunay{X} \to \Rspace$ the function that maps
each simplex to the radius of this circumball
and observe that $\Rf (P) \leq \Rf (Q)$ whenever $P \subseteq Q$. 
With probability one, $\Delaunay{X}$ is a simplicial complex
and $\Rf$ is a \emph{generalized discrete Morse function}.
In other words, there is a partition of $\Delaunay{X}$ into {intervals},
 where an \emph{interval} is a maximal set of simplices with common lower
bound and upper bound, $[L,U] = \{Q \mid L\subseteq Q \subseteq U\}$,
such that $\Rf(P)=\Rf(Q)$ for any two simplices in the interval,
and there is no further simplex contained in $L$ or containing $U$
with the same radius.
We call $[L,U]$ \emph{singular} if $L=U$ and \emph{non-singular} if $L$
is a proper subset of $U$.
In the former case, the unique simplex in the singular interval
is referred to as a \emph{critical simplex} of $\Rf$.
Given integers $\ell \leq k$ and a radius $r$,
we are interested in the number of intervals $[L,U]$ with
$\dime{L} = \ell$, $\dime{U} = k$, $\Rf (L) \leq r$,
and center in a Borel set $\Region \subseteq \Rspace^n$
of $n$-dimensional volume $\norm{\Region}$.
Here, the \emph{center} of an interval is the center of the circumball
of its simplices.

Our main result is given in terms of a certain expected volume of a random simplex.
Denote by $\uuu=(u_0, u_1, \ldots, u_k)$ a sequence of $k+1$ random points
chosen according to the uniform distribution on the unit sphere in $\Rspace^k$,
and write $\Volume{\uuu}$ for the $k$-dimensional volume of the $k$-simplex
spanned by the $u_i$.
The essential ingredient in the following is the expected volume
$\Expected{ \Volume{\uuu}^{n-k+1} \One{k-\ell} (\uuu) }$, in which
\begin{align}
  \One{k-\ell} (\uuu)  &=  \left\{ \begin{array}{ll}
                             1  &  \mbox{\rm if~} k-\ell
                                   \mbox{\rm ~of the~} k+1
                                   \mbox{\rm ~facets are visible from~} 0 , \\
                             0  &  \mbox{\rm otherwise} .
                           \end{array} \right.
\end{align}

\begin{theorem}[Main Result]
  \label{thm:main}
  Let $X$ be sampled from a Poisson point process with density $\density > 0$
  in $\Rspace^n$.
  For any $r > 0$, the expected number of intervals of the Poisson--Delaunay
  mosaic with center in a Borel set $\Region$ and radius at most $r$
  is given by the lower incomplete Gamma function,
  \begin{align}
    \frac{\gamma (k, \density \nu_n r^n)}{\Gamma (k)}
              \Ccon{\ell}{k}{n} \cdot \density \norm{\Region} ,
  \end{align}
  in which $\nu_n$ is the volume of the $n$-dimensional unit ball,
  $\sigma_n = n \nu_n$ is the $(n-1)$-dimensional surface area of the unit sphere,
  and the constant is explicitly given by
  \begin{align}
    \Ccon{\ell}{k}{n}  &=
      \tfrac{\sigma_n \cdot \sigma_{n-1} \cdot \ldots \cdot \sigma_{n-k+1}}
            {\sigma_1 \cdot \sigma_2 \cdot \ldots \cdot \sigma_k}
      \tfrac{\Gamma(k) n^{k-1} k!^{n-k} \sigma_k^{k+1}}
            {(k+1) \sigma_n^k}
        \; \Expected{ \Volume{\uuu}^{n-k+1} \One{k-\ell} (\uuu) } .
  \end{align}
  For dimension $n=2,3,4$, the expectation can be computed explicitly,
  and the resulting numerical values for $\Ccon{\ell}kn$
  are given in Table \ref{tbl:constants}.
\end{theorem}
\begin{table}[hbt]
  \centering
  \small{ \begin{tabular}{r||rrr|rrrr|rrrrr}
    \multicolumn{1}{c||}{$\Ccon{\ell}{k}{n}$} & \multicolumn{3}{c|}{$n=2$}
            & \multicolumn{4}{c|}{$n=3$} & \multicolumn{5}{c}{$n=4$}  \\
      &  $k = 0$  &  $1$  &  $2$
      &      $0$  &  $1$  &  $2$  &  $3$       
      &      $0$  &  $1$  &  $2$  &  $3$  &  $4$              \\ \hline \hline
    $\ell = 0$      & $1.00$ & $0.00$ & $0.00$
                    & $1.00$ & $0.00$ & $0.00$ & $0.00$ 
                    & $1.00$ & $0.00$ & $0.00$ & $0.00$ & $0.00$ \\
           $1$      &        & $2.00$ & $1.00$ 
                    &        & $4.00$ & $2.55$ & $1.21$ 
                    &        & $8.00$ & $5.66$ & $3.55$ & $1.66$ \\
           $2$      &        &        & $1.00$ 
                    &        &        & $4.85$ & $3.70$
                    &        &        &$17.66$ &$18.96$ &$11.14$ \\
           $3$      &        &        &       
                    &        &        &        & $1.85$
                    &        &        &        &$15.40$ &$14.22$ \\
           $4$      &        &        &       
                    &        &        &        &       
                    &        &        &        &        & $4.74$
  \end{tabular} }
  \caption{Approximations of the constants in the expected numbers
    of critical simplices (\emph{diagonal})
    and non-singular intervals (\emph{off-diagonal}) for a Poisson point process
    in $\Rspace^2$ on the \emph{left}, in $\Rspace^3$ in the \emph{middle},
    and in $\Rspace^4$ on the \emph{right}.}
  \label{tbl:constants}
\end{table}

Recall that the \emph{lower incomplete Gamma function} is defined as
$\gamma (k, x)  =  \int_{t=0}^x t^{k-1} e^{-t} \diff t$.
Its better known complete version, $\Gamma (k) = \gamma (k, \infty)$,
evaluates to $(k-1)!$ for every positive integer $k$.
For half-integers, we have
$\Gamma (k+\tfrac{1}{2}) = \tfrac{(2k)!}{4^k k!} \sqrt{\pi}$.
The $n$-dimensional volume of the unit ball in $\Rspace^n$ is
$\nu_n=\sqrt{\pi}^n / \Gamma (\tfrac{n}{2} + 1)$,
and the $(n-1)$-dimensional volume of its boundary
is $\sigma_n=2\sqrt{\pi}^n / \Gamma (\tfrac{n}{2})$.

\begin{figure}[hbt]
  \centering
  {\includegraphics[width=0.3\textwidth]{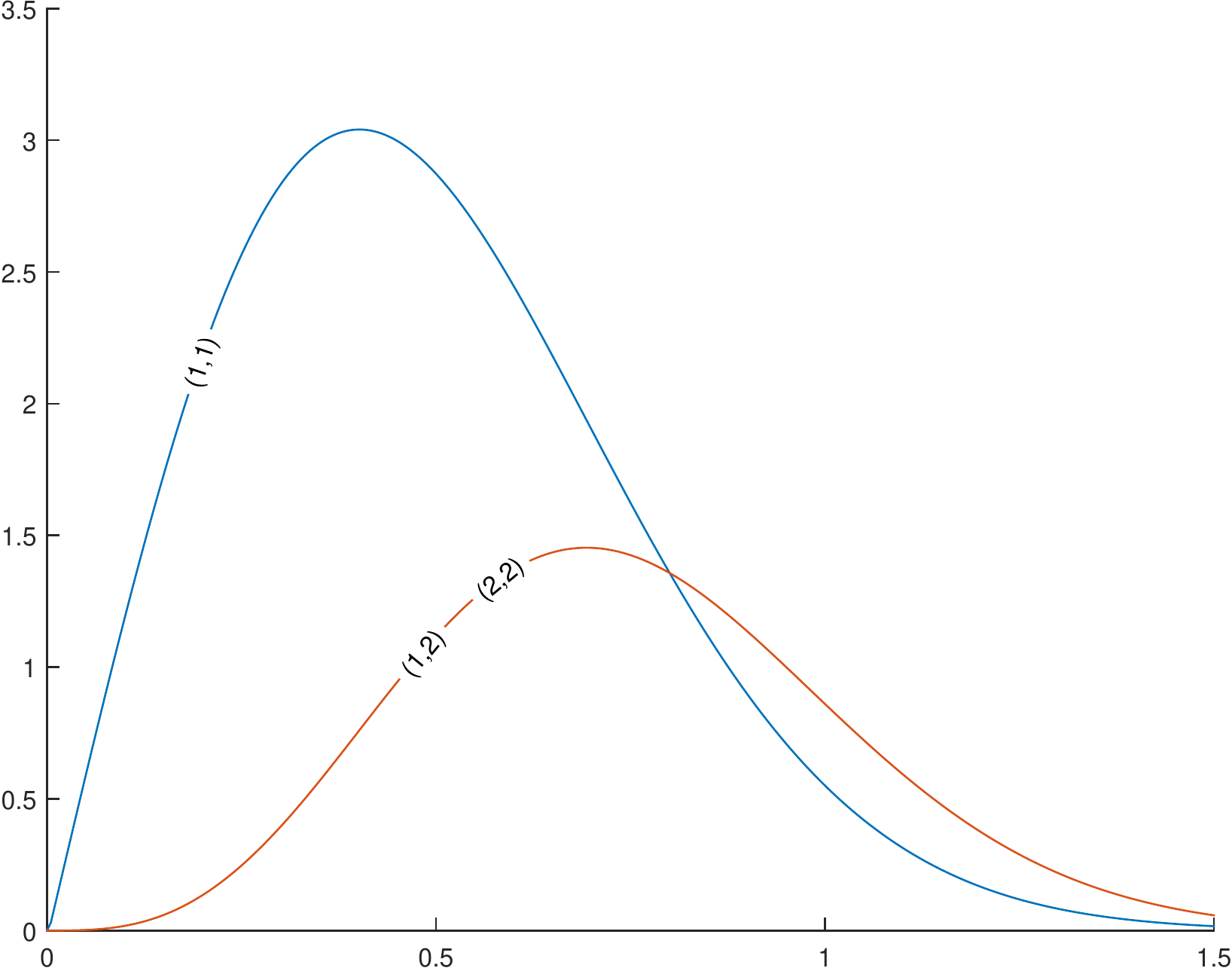}
   \includegraphics[width=0.3\textwidth]{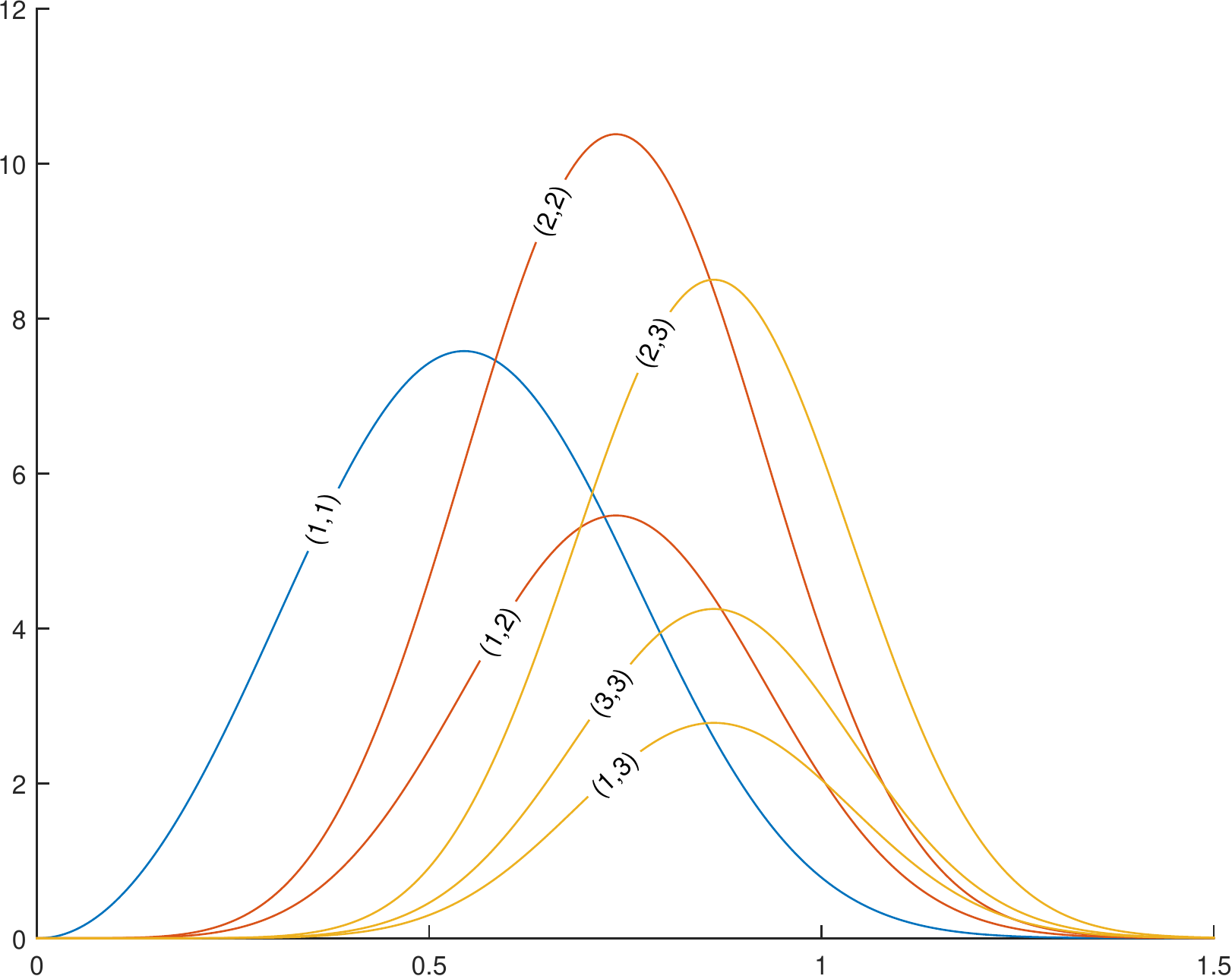}
   \includegraphics[width=0.3\textwidth]{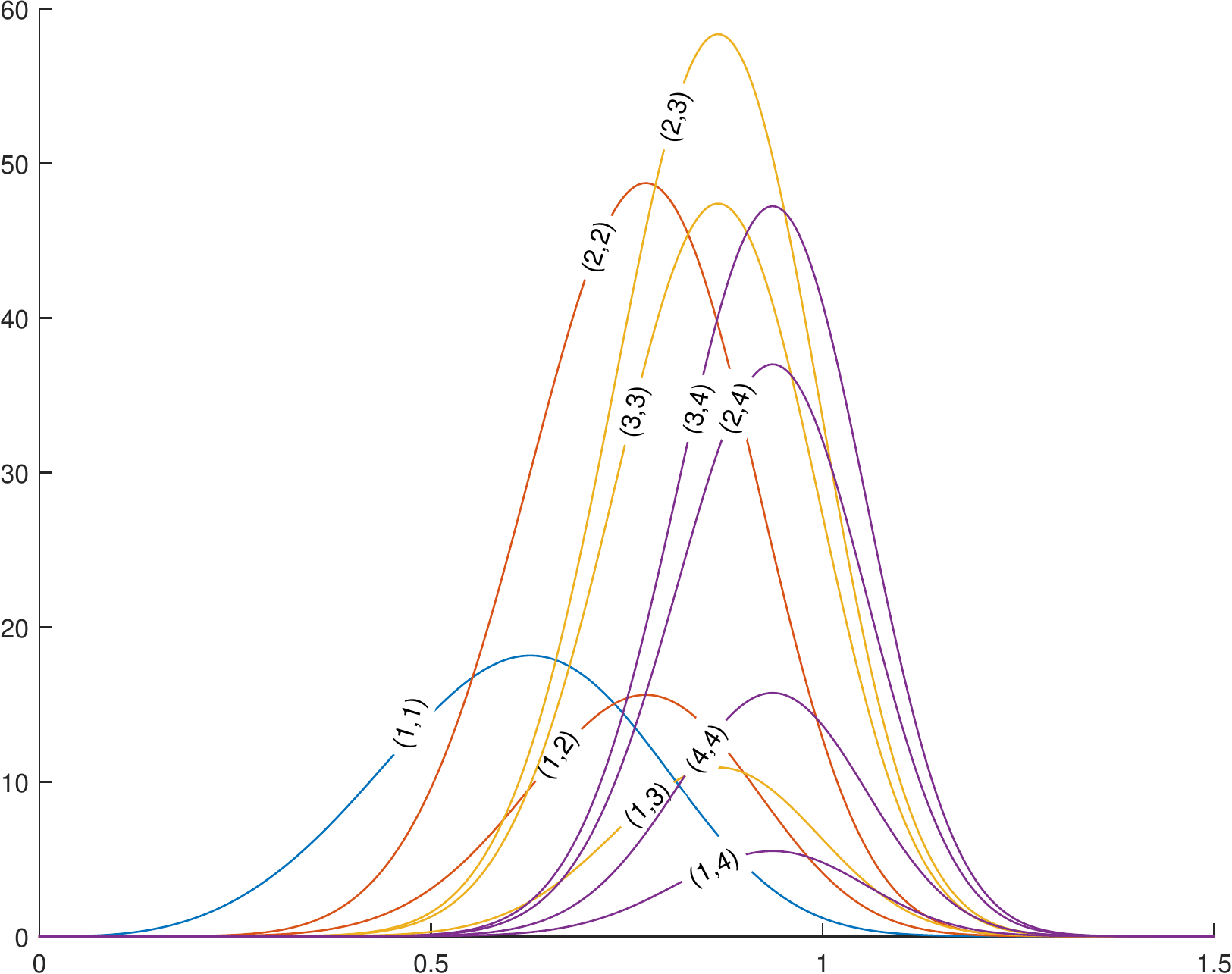}}
  {\includegraphics[width=0.3\textwidth]{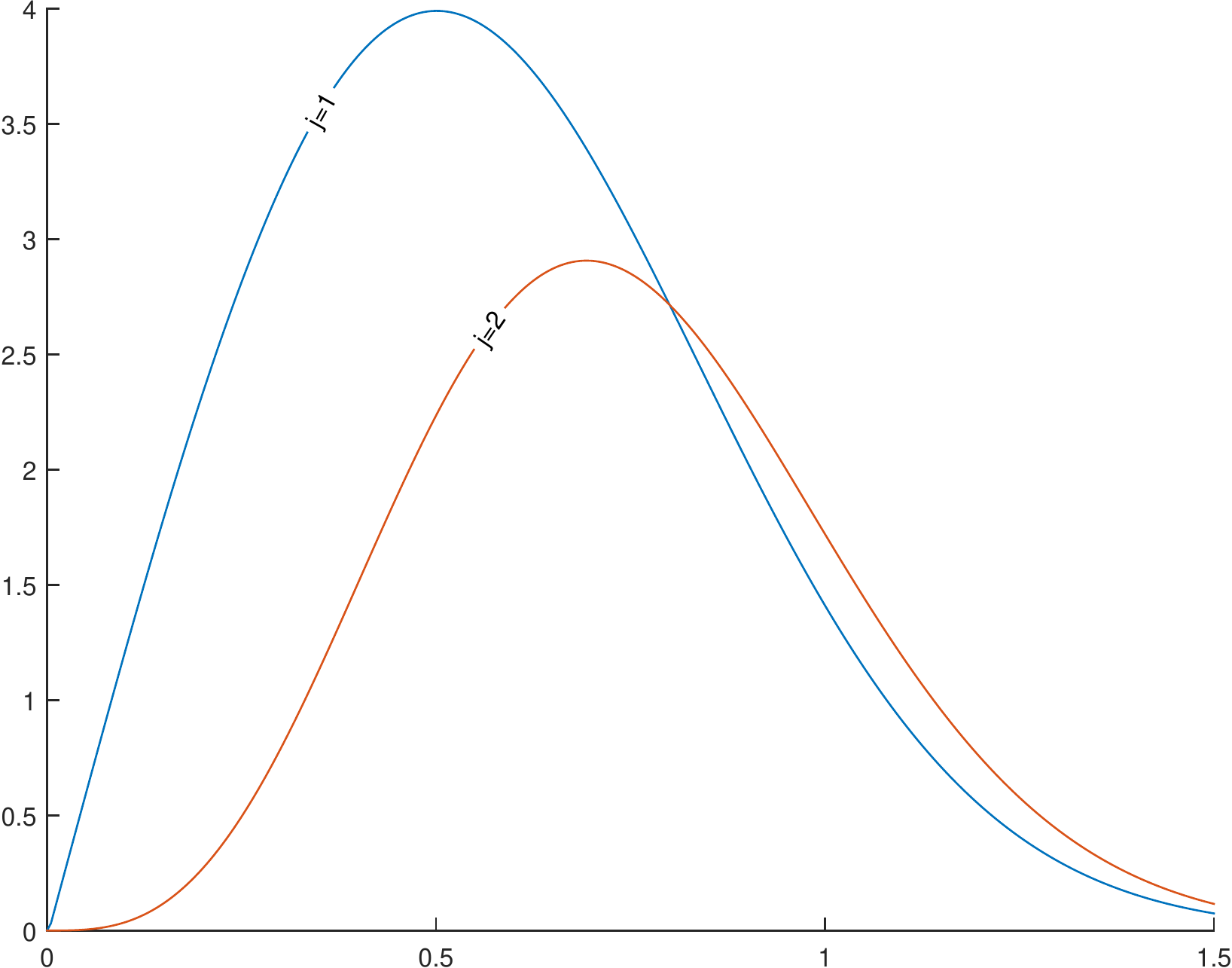}
   \includegraphics[width=0.3\textwidth]{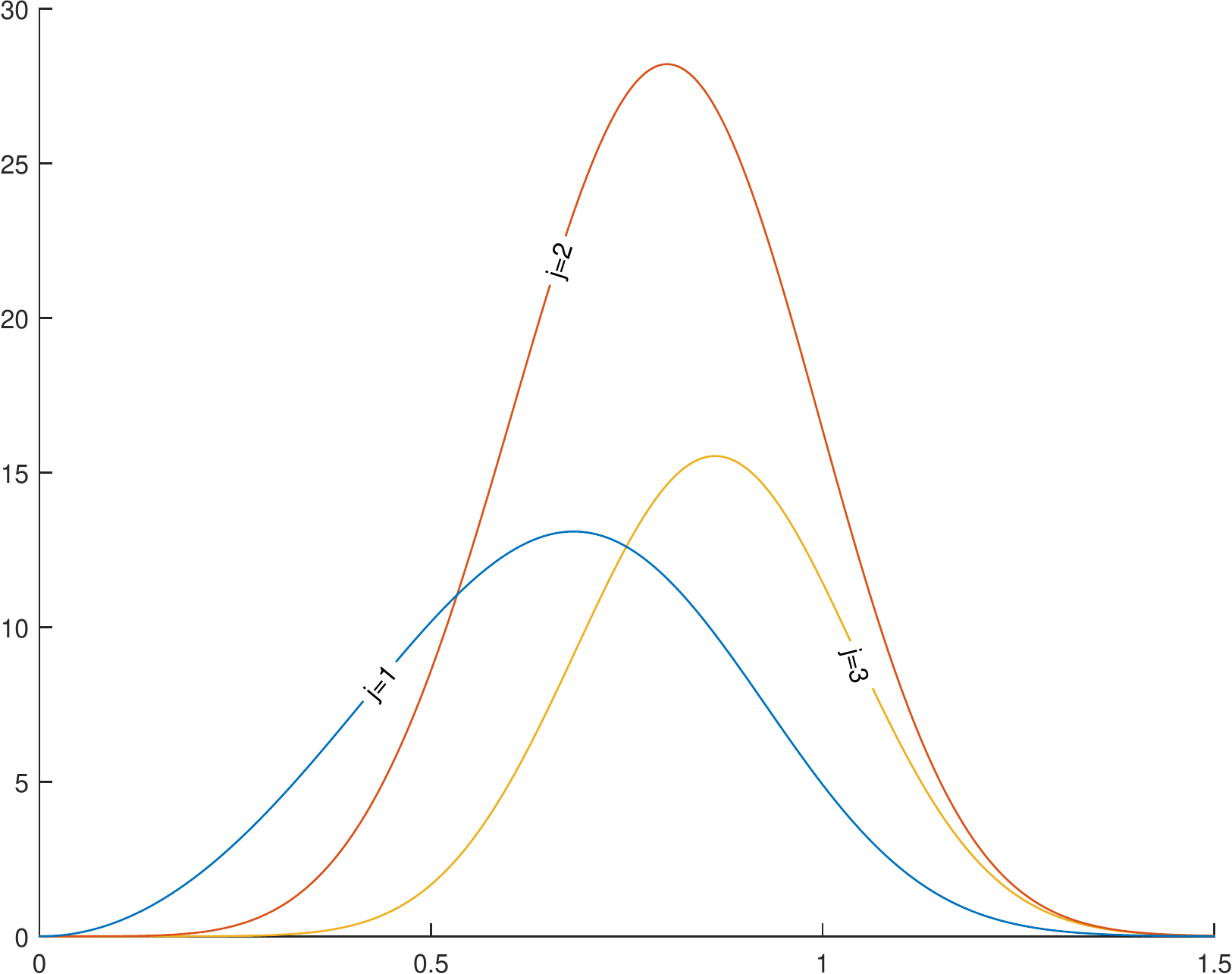}
   \includegraphics[width=0.3\textwidth]{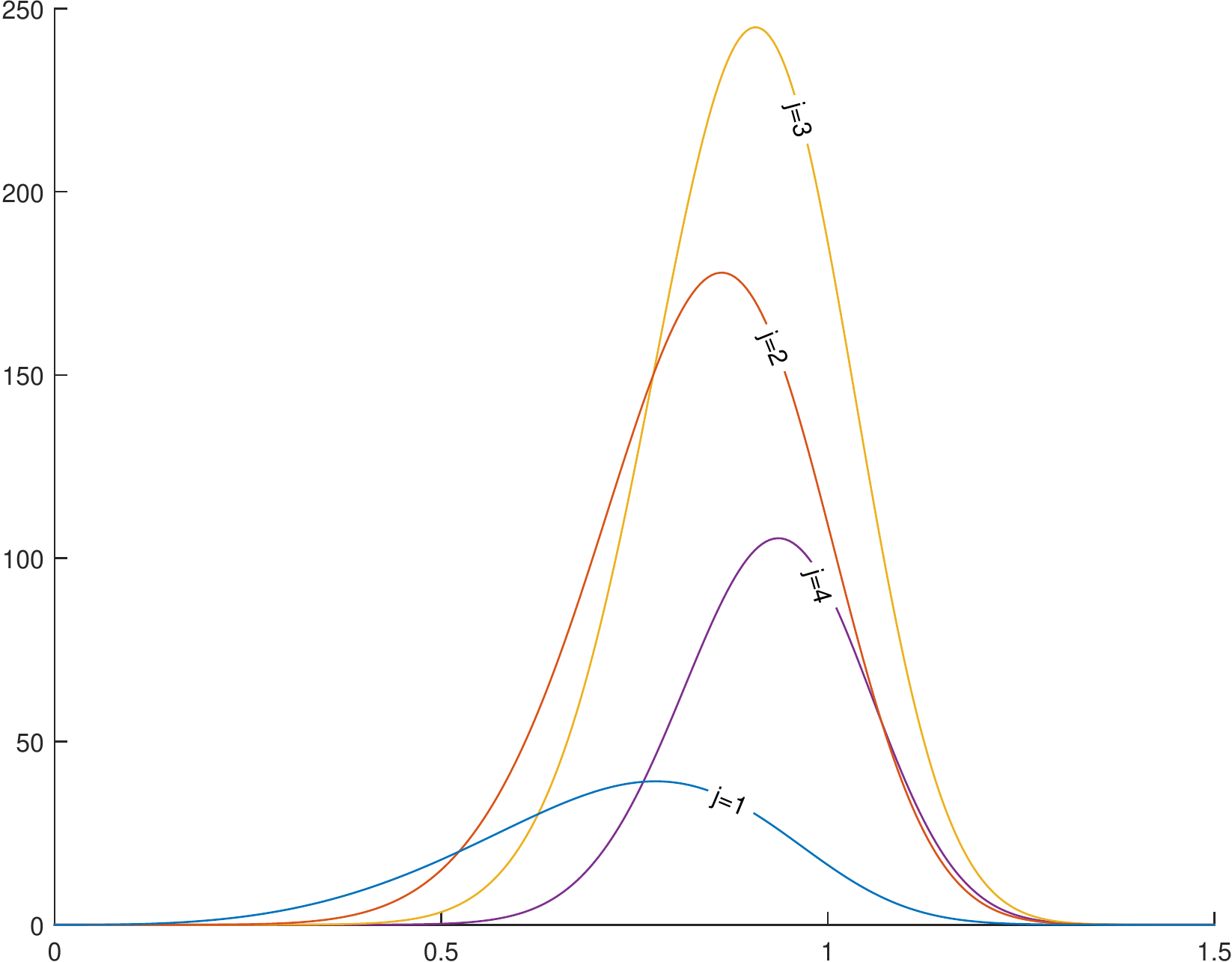}}
  \caption{\emph{Top}: the distribution of the expected number of intervals 
    as a function of the radius.
    The graphs are obtained by drawing $\Ccon{\ell}{k}{n}$ times
    the derivative of $\gamma (k, \nu_n r^n)$ normalized by $\Gamma (k)$,
    for $1 \leq \ell \leq k \leq n$,
    with $n = 2, 3, 4$ from \emph{left} to \emph{right}.
    \emph{Bottom}:  the corresponding distribution of the expected number of
    Delaunay simplices.
    We assume $\density = 1$ in all cases.}
  \label{fig:gammas}
\end{figure}

Let us make Theorem \ref{thm:main} more precise by introducing some notation.
The centers of the circumballs of the simplices form a point process
in $\Rspace^n$, which in general is not simple.
Yet if we restrict this to the centers of intervals --- and thus merge
the centers of the circumballs of all its simplices into one ---
we obtain a simple point process.
The expected number of centers of intervals in a region $\Region$
is the \emph{intensity measure} of the interval process.
Since the underlying Poisson point process is translation invariant,
so is the intensity measure. 
This shows that the expected number of intervals with 
$\dime{L} = \ell$, $\dime{U} = k$, $\Rf (L) \leq r$,
and center in $\Region $ factorizes into
$G^n_{\ell,k}(r) \Ccon{\ell}{k}{n} \cdot \rho \norm{\Region}$,
with some constant $\Ccon{\ell}{k}{n}$ called the \emph{intensity}
of the interval process, and where $G^n_{\ell,k}$ is an increasing function
with $\lim_{r \to \infty} G^n_{\ell,k}(r) = 1$.
Hence, $G^n_{\ell,k}$ is a probability distribution function,
and we call $G^n_{\ell,k}$ the distribution of the radius
of a \emph{typical interval}.
Theorem \ref{thm:main} shows that the radius of the typical interval
is Gamma distributed.
Clearly, the expected total number of intervals
in the Poisson--Delaunay mosaic with center in $\Region$ is
$\sum_{k=0}^n \sum_{\ell = 0}^k \Ccon{\ell}{k}{n} \cdot \density \norm{\Region}$.
More interestingly, the centers of the circumballs of the $j$-dimensional
simplices in the Poisson--Delaunay mosaic in $\Rspace^n$ is again a point process.
Its intensity is
\begin{align}
  \Dcon{j}{n} &= \sum_{k=j}^n \sum_{\ell=0}^j \binom{k-\ell}{k-j} \Ccon{\ell}{k}{n},
\end{align}
which can be evaluated explicitly for $n=2,3,4$; see Table \ref{tbl:PD-constants}.
This extends the result of Miles mentioned in \cite{ScWe08} to $n=4$.
By Theorem \ref{thm:main}, we can also get the expected number
of simplices with circumradius at most $r$ and center in $\Region$.
This seems to be one of the rare examples in the theory of random complexes,
in which the precise distribution can be computed.
\begin{corollary}[Delaunay Simplices]
  \label{Cor:dist-simplex}
  Let $X$ be sampled from a Poisson point process with density $\density > 0$
  in $\Rspace^n$.   
  The expected number of $j$-dimensional simplices of
  the Poisson--Delaunay mosaic with circumradius at most $r$
  and center in $\Region$ is 
  \begin{align}
    G^n_{j} (r) \Dcon{j}{n} \cdot \density \norm{\Region}
      &=  \sum_{k=j}^n \frac{\gamma (k, \density \nu_n r^n)}{\Gamma (k)}
          \sum_{\ell = 0}^j \binom{k-\ell}{k-j}
                            \Ccon{\ell}{k}{n} \cdot \density \norm{\Region} .
     \label{eq:numbsimpl}
  \end{align}
  Hence, the distribution of the circumradius of the typical
  $j$-dimensional simplex is a mixed Gamma distribution:
  \begin{align}
    G^n_{j} (r)  &=  \sum_{k=j}^n \frac{\gamma (k, \density \nu_n r^n)}{\Gamma (k)} 
                     \sum_{\ell = 0}^j \binom{k-\ell}{k-j}
                                       \frac{\Ccon{\ell}{k}{n}}{\Dcon{j}{n}} .
    \label{eq:typ-j-simpl}
  \end{align}
  For dimension $n=2,3,4$, the constants can be computed explicitly,
  with values given in Table \ref{tbl:PD-constants};
  see also Figure \ref{fig:gammas}. 
\end{corollary}

\begin{table}[hbt]
  \centering
  \small{ \begin{tabular}{r||rrrrr}
    \multicolumn{1}{c||}{$\Dcon{j}{n}$}
             & $j=0$  & $1$     & $2$    & $3$    & $4$       \\ \hline \hline
    $n = 2$  & $1.00$ & $3.00$  & $2.00$ &        &           \\
        $3$  & $1.00$ & $7.76$  &$13.53$ & $6.76$ &           \\
        $4$  & $1.00$ & $18.88$ &$65.55$ &$79.44$ &$31.77$
  \end{tabular}}
  \caption{Approximations of the constants in the expected numbers
    of simplices in a Poisson--Delaunay mosaic.
    They are straightforward in two dimensions,
    they have been found by R.\ Miles \cite{ScWe08} in three dimensions,
    and except for $j = 0, 3, 4$ they are new in four dimensions.}
  \label{tbl:PD-constants}
\end{table}

As an application we obtain that the expected number of $j$-dimensional
simplices in the Delaunay complex for radius $r$ with center in $\Region$
is given by \eqref{eq:numbsimpl}.
For $j=n$, this is a consequence of the Complementary Theorem due to Miles
\cite{Mil74}, and we mention the close relation to Gamma type results
by Baumstark and Last \cite{BaLa09} and Chenavier \cite{Che14}.
It should also be pointed out that Corollary \ref{Cor:dist-simplex}
can be converted to results for the dual tessellation,
the Poisson--Voronoi tessellation.
Then \eqref{eq:numbsimpl} gives the intensity of the $(n-j)$-dimensional
face process of the Poisson--Voronoi tessellation,
and \eqref{eq:typ-j-simpl} the distance of the typical $(d-j)$-face
to the closest point of the Poisson point process.

The main technical tool used to prove Theorem \ref{thm:main} is a
Blaschke--Petkantschin type formula that integrates over spheres.
Let $0 < k \leq n$, write $\LGrass{k}{n}$ for the Grassmannian
consisting of all $k$-planes passing through the origin in $\Rspace^n$,
let $S = S(L)$ be the $(k-1)$-dimensional unit sphere
inside $L \in \LGrass{k}{n}$,
let $\uuu = (u_0, u_1, \ldots, u_k)$ be a sequence of $k+1$ points,
and write $\Volume{\uuu}$ for the $k$-dimensional volume of
the $k$-simplex spanned by the $u_i$.
Then for every non-negative function $f$ on $k+1$ points, we have
\begin{align}
  \int\displaylimits_{\xxx \in (\Rspace^n)^{k+1}} \!\!\! f(\xxx) \diff \xxx  &=
    \int\displaylimits_{L \in \LGrass{k}{n}}
    \int\displaylimits_{z \in \Rspace^n}
    \int\displaylimits_{r \geq 0}
    \int\displaylimits_{\uuu \in S^{k+1}} \!\!\!
      f(z + r \uuu) r^{nk-1} (k! \Volume{\uuu})^{n-k+1}
                    \diff \uuu \diff r \diff z \diff L .
  \label{eqn:BP}
\end{align}
This equation generalizes Theorem 7.3.1 in \cite[page 287]{ScWe08}
to $k \leq n$.

While Theorem \ref{thm:main} makes a statement
about expectations in a fixed region $\Omega$,
a standard ergodic argument implies that for a sequence of regions
$\Region_1 \subseteq \Region_2 \subseteq \ldots$ covering the entire space,
the numbers of intervals inside $\Region_i$, normalized by $\norm{\Region_i}$,
converge to the corresponding constants almost surely as random variables, see 
\cite{Mil70} for details.

\ourparagraph{Outline.}
Section \ref{sec:2} introduces the background on discrete Morse
theory and Poisson--Delaunay mosaics.
Section \ref{sec:3} explains the essential integral geometric tools
used to prove Theorem \ref{thm:main}.
Section \ref{sec:4} presents combinatorial and probabilistic results on
inscribed simplices.
Section \ref{sec:5} uses these results to compute the constants
given in Tables \ref{tbl:constants} and \ref{tbl:PD-constants}.
Section \ref{sec:6} concludes the paper.

\section{Poisson--Delaunay Complexes}
\label{sec:2}

In this section, we introduce the necessary background on Poisson point
processes and the radius function on Delaunay mosaics.

\ourparagraph{Poisson point processes.}
We study properties of randomly generated discrete point sets in $\Rspace^n$
using a \emph{Poisson point process} with \emph{density} $\density > 0$,
which can be characterized by the following two properties:

\begin{enumerate}\denselist
  \item[1.] The numbers of points sampled within a finite collection
    of pairwise disjoint Borel sets are independent random variables;
  \item[2.] The expected number of points sampled within a Borel
    set is $\density$ times the Lebesgue measure of the set.
\end{enumerate}

\noindent
See \cite{Kin93} for a good introduction to Poisson point processes.
Writing $X \subseteq \Rspace^n$ for the set sampled from the process,
we can express Condition 2 more succinctly as
$\Expected{\card{X \cap B}} = \density \norm{B}$.
The two conditions imply that the number of points sampled in
a Borel set $B$ has a Poisson distribution with parameter $\density \norm{B}$.
In particular,
the probability of having no point in $B$ is
$\Probable{}{X \cap B = \emptyset}  =  e^{- \density \norm{B}}$.
Another property that will be important in this paper is the following.
\begin{lemma}[General Position]
  \label{lem:genpos}
  Let $X$ be sampled from a Poisson point process in $\Rspace^n$.
  With probability $1$, $X$ is a countable set of points such that
  for every $0 \leq k < n$,
  \begin{enumerate}\denselist
    \item[1.] no $k+2$ points belong to a common $k$-plane,
    \item[2.] no $k+3$ points belong to a common $k$-sphere,
    \item[3.] considering the unique $k$-sphere that passes through
      $k+2$ points, no $k+1$ of these points belong to a $k$-plane
      that passes through the center of the $k$-sphere.
  \end{enumerate}
\end{lemma}

\noindent
It is not difficult to prove the above lemma,
and we refer to \cite{Mil70} and \cite{ScWe08} for further information.
We say that $X$ is in \emph{general position} if it satisfies the
conditions stated in Lemma \ref{lem:genpos}.
Since it happens with probability $1$,
we will assume that $X$ is in general position throughout this paper.

\ourparagraph{Radius function on Delaunay mosaic.}
The \emph{Voronoi domain} of a point $x \in X$ consists of all points
for which $x$ minimizes the Euclidean distance:
$\Voronoi{x}  =  \{ a \in \Rspace^n  \mid  \Edist{a}{x} \leq \Edist{a}{y},
          \mbox{\rm ~for all~} y \in X \}$.
With probability $1$, every Voronoi domain is a bounded
convex polyhedron \cite{ScWe08}.
The \emph{Voronoi diagram} is the collection of Voronoi domains,
and the \emph{Delaunay mosaic} is the nerve of the Voronoi domains;
that is:
the collection of simplices $Q \subseteq X$ with non-empty common
intersection of Voronoi domains, $\Voronoi{Q} = \bigcap_{x \in Q} \Voronoi{x}$.
The Delaunay mosaic is a simplicial complex iff the Voronoi diagram
is \emph{primitive}, that is: the intersection of any
$1 \leq k+1 \leq n+2$ Voronoi domains is either empty or $(n-k)$-dimensional.
In particular, the intersection of any $n+2$ domains is necessarily empty.
Since $X$ is in general position with probability $1$,
the Voronoi diagram is primitive with probability $1$,
and the Delaunay mosaic is a simplicial complex,
again with probability $1$.

By construction, every point $z \in \Voronoi{Q}$ is equally far from
all points in $Q$ and at least as far from all points in $X \setminus Q$.
We call the sphere with center $z$ and radius $\Edist{z}{x}$, $x \in Q$,
a \emph{circumsphere} because all points of $Q$ lie on the sphere,
and we call it \emph{empty} because all points of $X$ lie on or outside the 
sphere.
The \emph{radius function}, $\Rf \colon \Delaunay{X} \to \Rspace$,
maps every simplex to the smallest radius of all its empty circumspheres.
As proved in \cite{BaEd15}, the radius function satisfies the conditions of
a generalized discrete Morse function.
We refer to \cite{For98} for an introduction to discrete Morse theory
and to \cite{Fre09} for the generalization needed here.
Aiming at a geometric characterization of the intervals,
we note that every simplex $Q \subseteq X$ has a unique
smallest circumsphere, namely the unique circumsphere whose center, $z$,
lies in the affine hull of $Q = \{x_0, x_1, \ldots, x_k\}$.
This circumsphere may or may not be empty.
Writing $z$ in terms of barycentric coordinates,
$z = \sum_{i=0}^k \zeta_i x_i$ with $\sum_{i=0}^k \zeta_i = 1$,
we note that the $\zeta_i$ are unique
and assuming general position they are also non-zero.
Interpreting the barycentric coordinates geometrically, we note that
$\zeta_i < 0$ iff the facet $Q \setminus \{x_i\}$ is \emph{visible}
from $z$; that is: there is a ray emanating from $z$ that first intersects
the interior of the facet before it intersects the interior of the $k$-simplex.
\begin{lemma}[Geometric Characterization of Intervals]
  \label{lem:intervals}
  Letting $X \subseteq \Rspace^n$ be in general position,
  a pair of simplices $L \subseteq U$ in $X$ defines an
  interval of $\Rf \colon \Delaunay{X} \to \Rspace$
  iff the smallest circumsphere of $U$ is empty
  and $L$ is the largest face of $U$ common to all facets that
  are visible from the center of the sphere.
\end{lemma}

\noindent
We refer to \cite{BaEd15} for a proof of this characterization.
Recall the special case of a critical simplex, $L = U$,
which is characterized by containing $z$ inside the convex hull.
In this case, the smallest circumsphere is also the smallest
sphere that encloses $U$.

\ourparagraph{Euler characteristic.}
Letting $K$ be a finite subset of $\Delaunay{X}$,
we write $\dcon{i}{} = \dcon{i}{} (K)$ for the number of $i$-simplices,
for $0 \leq i \leq n$.
The sum of these numbers is the \emph{size} and the
alternating sum is the \emph{Euler characteristic}:
$\card{K}  =  \sum_{i=0}^n \dcon{i}{} (K)$ and
$\chi (K)  =  \sum_{i=0}^n (-1)^i \dcon{i}{} (K)$.
Importantly, if $K$ is a complex, then the Euler characteristic
is determined by its homotopy type.
For example, if $K$ is contractible, then $\chi (K) = 1$.
Now suppose that $K$ is a union of intervals of the radius function,
and let $\ccon{i}{}{} = \ccon{i}{}{} (K)$ be the number of $i$-dimensional
critical simplices of $\Rf$ restricted to $K$.
The discrete version of the Morse Relation asserts that the
Euler characteristic is also the alternating sum of the $\ccon{i}{}{}$:
\begin{lemma}[Discrete Morse Relation]
  \label{lem:discreteMorse}
  Letting $X \subseteq \Rspace^n$ be in general position, and
  $K$ a finite subset of the Delaunay mosaic that is a union
  of intervals of the radius function,
  the Euler characteristic satisfies
  $\chi (K) = \sum_{i=0}^n (-1)^i \ccon{i}{}{} (K)$.
\end{lemma}
\ourproof
  Consider an interval $[L,U]$, and assume that
  $k = \dime{U} > \ell = \dime{L}$.
  For each $0 \leq j \leq k - \ell$, the number of $(\ell+j)$-simplices
  in the interval is the number of subsets of size $j$ of $U \setminus L$.
  The contribution of the interval to the Euler characteristic is
  therefore $\pm \sum_{j=0}^{k-\ell} (-1)^j \binom{k-\ell}{j} = 0$.
  It follows that the non-singular intervals have no effect on the
  Euler characteristic, which implies the claimed relation.
\eop

\begin{figure}[hbt]
  \centering \resizebox{!}{2.0in}{\includegraphics{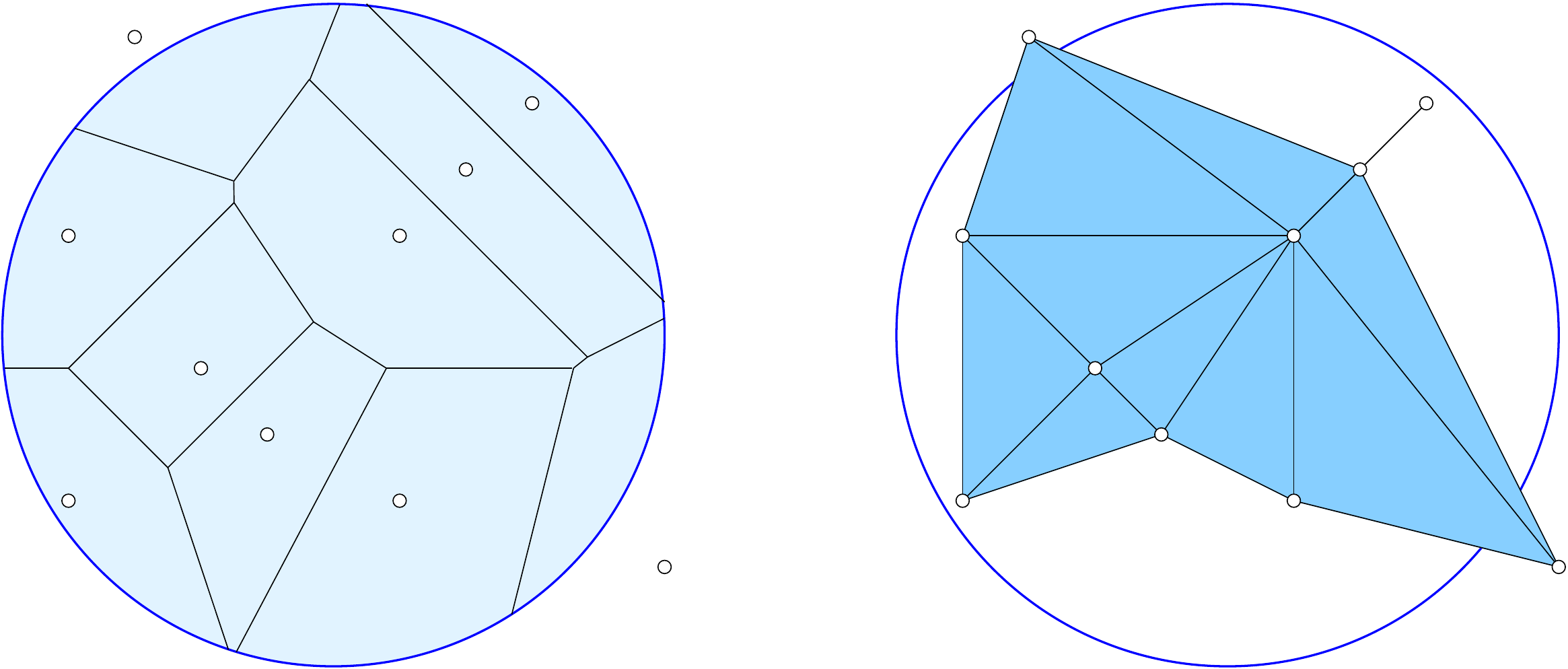}}
  \caption{The Voronoi diagram restricted to a disk on the \emph{left},
    and the corresponding restricted Delaunay mosaic, $K_0$, on the 
\emph{right}.
    In this case, the set $K_1$ consists of all simplices in $K_0$
    except for the two vertices that lie outside $\Region$,}
  \label{fig:complexes}
\end{figure}

\ourparagraph{Subsets and subcomplexes.}
Assuming $X$ is in general position, 
we use $\Region$ to specify two subsets of the Delaunay mosaic.
\begin{itemize}\denselist
  \item The subcomplex $K_0 = K_0 (R)$ of $\Delaunay{X}$ consists
    of all simplices $Q$ such that $\Voronoi{Q} \cap \Region \neq \emptyset$;
    see Figure \ref{fig:complexes}.
    Equivalently, $K_0$ consists of all simplices such that
    $\bigcap_{x \in Q} [\Voronoi{x} \cap \Region] \neq \emptyset$.
    Since the $\Voronoi{x} \cap \Region$ are convex,
    the Nerve Theorem applies and asserts that $K_0$ and $\Region$
    have the same homotopy type; see e.g.\ \cite[page 59]{EdHa10}.
    Since $\Region$ is contractible, this implies $\chi (K_0) = 1$.
  \item The subset $K_1 = K_1 (R)$ of $\Delaunay{X}$ consists
    of all simplices in $K_0$ whose smallest empty circumspheres
    have the center in $\Region$.
    We can construct $K_1$ by removing one simplex at a time from $K_0$.
    Each removed simplex changes the Euler characteristic by $1$,
    which gives 
    \begin{align}
      \label{eqn:Euler-D}
      1 - \card{K_0 \setminus K_1}  \leq  &\chi (K_1)
                                    \leq  1 + \card{K_0 \setminus K_1} .
    \end{align}
    We will work with $K_1$ throughout this paper, in particular Theorem \ref{thm:main}
    is proved for $K_1$.
\end{itemize}
$K_0$ is interesting for nice sets: in the appendix we prove that if $\Omega$ is
a ball, then the difference between the number of simplices in $K_0$ and $K_1$
is $\littleoh{\norm{\Region}}$. This is so because the simplices in
$K_0 \setminus K_1$ correspond to common intersections of Voronoi domains that
touch the boundary of $\Region$, and there are not too many of them.
Together with \eqref{eqn:Euler-D}, this property implies
$\chi (K_1) = \littleoh{\norm{\Region}}$;
see Appendix \ref{app:A} for details.
The proof straightforwardly extends from balls to any region with $C^2$-boundary.
It implies that Theorem \ref{thm:main} also holds for $K_0$ for smooth $\Region$
after adding $\littleoh{1} \density \norm{\Region}$.

Theorem \ref{thm:main} has the same constant for any Borel set $\Omega$, so
it suffices to obtain them for $\Omega$ a ball.
We can therefore use the aforementioned relation,
$\chi (K_1) = \littleoh{\norm{\Region}}$, which implies the Euler
relation for the constants: $\sum_{j = 0}^{n} (-1)^j \Dcon{j}{n} = 0$. 
Similarly, we have $\sum_{j = 0}^{n} (-1)^j \Ccon{j}{j}{n} = 0$
because of Lemma \ref{lem:discreteMorse}.
Later we will also use
the fact that $\Dcon{n-1}{n} = \tfrac{n+1}{2} \Dcon{n}{n}$, which is clear
for $K_0$ for balls and hence also always holds for $K_1$.

\section{Integral Geometry}
\label{sec:3}

We obtain all our results by combining relations from discrete geometry
with equations from integral geometry.
This section introduces the latter.

\ourparagraph{Spherical Blaschke--Petkantschin Formula.}
In its basic form, the \emph{Bla\-schke--Pet\-kan\-tschin Formula} decomposes
an integral over $\Rspace^n$ into an integral over the
Grassmannian, $\LGrass{k}{n}$,
times an integral over a $k$-plane and its orthogonal $(n-k)$-plane.
We start with the form  given in \cite[Equation (27)]{Mil71}:
\begin{align}
  \int\displaylimits_{\xxx \in (\Rspace^n)^{k+1}} f(\xxx) \diff \xxx
     &=  \int\displaylimits_{L \in \LGrass{k}{n}}
         \int\displaylimits_{h \in L^\perp}
         \int\displaylimits_{\xxx \in L^{k+1}} 
         f(h+\xxx) (k! \Volume{\xxx})^{n-k} \diff \xxx \diff h \diff L ,
  \label{eqn:BP1}
\end{align}
in which $\xxx = (x_0, x_1, \ldots, x_k)$, each $x_i$ is a point in $\Rspace^n$,
$\Volume{\xxx}$ is the $k$-dimensional volume of the simplex spanned by $\xxx$,
and $f$ is a non-negative function on $k+1$ points.
The $(n-k)$-th power of the volume
compensates for the biased measure on $\Rspace^n$ introduced by
the Grassmannian.
Using Theorem 7.3.1 in \cite[page 287]{ScWe08},
we expand the innermost integral into
\begin{align}
  k! ~\int\displaylimits_{z \in L} ~~
     \int\displaylimits_{r \geq 0}
     \int\displaylimits_{\uuu \in S^{k+1}}
       r^{k^2-1} \Volume{\uuu} f(h+z+r\uuu)
       (k! \Volume{z+r\uuu})^{n-k} \diff \uuu \diff r \diff z ,
  \label{eqn:BP2}
\end{align}
in which $S = S(L)$ is the unit $(k-1)$-sphere in $L$,
and $\uuu = (u_0, u_1, \ldots, u_k)$, with each $u_i$ a point on $S$.
Note that $\Volume{z+r\uuu} = r^k \Volume{\uuu}$, so we get
$k^2-1 + (n-k)k = nk-1$ as the final power of the radius.
Plugging \eqref{eqn:BP2} into \eqref{eqn:BP1} and joining
the integration over $L^\perp$ and $L$, we get
\begin{align}
  \int\displaylimits_{\xxx \in (\Rspace^n)^{k+1}} \!\! f(\xxx) \diff \xxx
     &=  \int\displaylimits_{L \in \LGrass{k}{n}}
         \int\displaylimits_{z \in \Rspace^n}
         \int\displaylimits_{r \geq 0}
         \int\displaylimits_{\uuu \in S^{k+1}}
           \!\! r^{nk-1} f(z+r\uuu) (k! \Volume{\uuu})^{n-k+1}
           \diff \uuu \diff r \diff z \diff L ,
  \label{eqn:BPagain}
\end{align}
which is Equation \eqref{eqn:BP} stated in the Introduction.
Not surprisingly, for $k = n$ this is the equation in
Theorem 7.3.1 of \cite{ScWe08}.

\ourparagraph{Slivnyak--Mecke Formula.}
In a nutshell, the \emph{Slivnyak--Mecke Formula} writes the expectation
of a random variable of a Poisson point process as an integral over
the space on which the process is defined;
see \cite[page 68]{ScWe08}.
Motivated by the characterization of intervals in Lemma \ref{lem:intervals},
we are interested in the expected number of $k$-simplices
with empty smallest circumsphere.
Writing $\Pempty{\xxx}$ for the probability that the smallest
circumsphere of the simplex spanned by
$\xxx = (x_0, x_1, \ldots, x_k) \in (\Rspace^n)^{k+1}$ is empty,
we get the expected number of $k$-simplices
with empty smallest circumspheres as $\tfrac{1}{(k+1)!} \density^{k+1}
    \int_{\xxx \in (\Rspace^n)^{k+1}} \Pempty{\xxx} \diff \xxx$.
This integral is, of course, infinite, but it can be made finite by introducing 
geometric constraints, which can be either understood as restrictions of the 
integration domain, or as the corresponding indicator functions as multipliers 
at $\Pempty{\xxx}$.
For example, to get the number of critical $k$-simplices,
we would multiply with $\One{0} (\xxx)$, which we set to $1$
if the center of the smallest circumsphere
lies inside the $k$-simplex spanned by $\xxx$, and to $0$ otherwise.
More generally, we define
\begin{align}
  \One{k-\ell} (\xxx)  &=  \left\{ \begin{array}{ll}
                             1  &  \mbox{\rm if~} k-\ell
                                   \mbox{\rm ~of the~} k+1
                                   \mbox{\rm ~facets are visible from~} 0 , \\
                             0  &  \mbox{\rm otherwise} .
                           \end{array} \right.
\end{align}
Similarly, we let $\One{\Region} (\xxx)$ be the indicator function
that the center of this sphere is in $\Region$,
and we let $\One{r_0} (\xxx)$ be the indicator function that
this sphere has radius at most $r_0$.
We then get the expected number of $k$-simplices
whose smallest circumsphere is empty,
whose radius is at most $r_0$,
and whose center lies in $\Region$ such $k - \ell$ facets are visible as
\begin{align}
  \Expected{\ccon{\ell}{k}{}, r_0}  &=  \tfrac{1}{(k+1)!} \density^{k+1}
     \int\displaylimits_{\xxx \in (\Rspace^n)^{k+1}} \Pempty{\xxx}
     \; \One{k-\ell} (\xxx) \; \One{\Region} (\xxx)
     \; \One{r_0} (\xxx) \diff \xxx .
  \label{eqn:SM}
\end{align}
We refer to \eqref{eqn:SM} as a \emph{corollary of the
Slivnyak--Mecke Formula}.
Using the Spherical Blaschke--Petkantschin Formula \eqref{eqn:BPagain},
we will rewrite this integral and absorb the latter two
indicator functions by limiting the domain.
The indicator functions that distinguish between different numbers
of visible facets will remain and require further attention.

\ourparagraph{Counting intervals.}
We combine the Spherical Blaschke--Petkantschin Formula \eqref{eqn:BPagain}
with the corollary of the Slivnyak--Mecke Formula \eqref{eqn:SM}
to get a general expression for the expected number of intervals.
Before following up, we recall that the measure of the Grassmannian is
$\norm{\LGrass{k}{n}}
  = \tfrac{\sigma_n \sigma_{n-1} \cdot \ldots \cdot \sigma_{n-k+1}}
          {\sigma_1 \sigma_2     \cdot \ldots \cdot \sigma_k}$,
in which $\sigma_i$ is the $(i-1)$-dimensional volume of the unit 
$(i-1)$-sphere.
For example, $\LGrass{n-1}{n}$ may be identified with the
set of normal directions and
$\norm{\LGrass{n-1}{n}} = \tfrac{\sigma_n}{2}$,
half the volume of the $(n-1)$-dimensional sphere.

The expression for the expected number of intervals will contain
the lower incomplete Gamma function, $\gamma (k, x)$, as a factor.
We prefer to explain how this comes about before stating the formula.
Write $A = r^n \nu_n$ for the $n$-dimensional volume of the ball with
radius $r$ in $\Rspace^n$, and note that
$\diff A = n r^{n-1} \nu_n \diff r$.
Hence,
\begin{align}
  (n \nu_n^k) r^{nk-1} \diff r 
    &=  (n \nu_n^k) (r^n)^{k-1} r^{n-1} \diff r
     =  (n \nu_n^k) \left( \tfrac{A}{\nu_n} \right)^{k-1}
                           \tfrac{1}{n \nu_n} \diff A
     =  A^{k-1} \diff A .
  \label{eqn:substitution}
\end{align}
We use \eqref{eqn:substitution} to substitute and get the integral
into the form of the Gamma function.
Specifically,
\begin{align}
  n \nu_n^k \int\displaylimits_{r=0}^{r_0}
                   \density^k e^{-\density A} r^{nk-1} \diff r  
     &=  \int\displaylimits_{\density A = 0}^{\density A_0}
                   (\density A)^{k-1} e^{- \density A} \diff (\density A) 
      =  \gamma (k, \density A_0) ,
  \label{eqn:gammasub}
\end{align}
in which $A_0 = r_0^n \nu_n$.
With this introduction, we are ready to state and prove the general formula
for the expected number of intervals.
\begin{lemma}[Counting Intervals]
  \label{lem:counting}
  Let $X$ be chosen from a Poisson point process with density $\density > 0$
  in $\Rspace^n$.
  For any $1 \leq \ell \leq k \leq n$ and $r_0 \geq 0$, the expected
  number of $k$-simplices whose smallest circumsphere is empty,
  whose radius is at most $r_0$, and whose center lies in $\Region$
  such that $k - \ell$ facets are visible from the center is
  \begin{align}
    \Expected{\ccon{\ell}{k}{}, r_0}
      &=  \density \norm{\Region} \tfrac{k!^{n-k}}{n (k+1) \nu_n^k}
          \norm{\LGrass{k}{n}} \cdot \gamma (k, \density \nu_n r_0^n)
          \int\displaylimits_{\uuu \in (\Sspace^{k-1})^{k+1}}
             \Volume{\uuu}^{n-k+1} \One{k-\ell} (\uuu) \diff \uuu .
    \label{eqn:counting}
  \end{align}
\end{lemma}
\ourproof
  We start by rewriting the corollary of the Slivnyak--Mecke Formula
  \eqref{eqn:SM} using
  the Spherical Blaschke--Petkantschin Formula \eqref{eqn:BPagain}:
  \begin{align}
    \Expected{\ccon{\ell}{k}{}, r_0}
      &=  \tfrac{\density^{k+1}}{(k+1)!}
          \int\displaylimits_{\xxx \in (\Rspace^n)^{k+1}} \Pempty{\xxx}
            \; \One{k-\ell}(\xxx) \; \One{\Region}(\xxx)
            \; \One{r_0}(\xxx) \diff \xxx \\
      &=  \tfrac{\density^{k+1}}{(k+1)!}
        \!\int\displaylimits_{L \in \LGrass{k}{n}}
          \int\displaylimits_{z \in \Region}
          \int\displaylimits_{r = 0}^{r_0}
          \int\displaylimits_{\uuu \in S^{k+1}}
          \!\!\!\!\!\! r^{nk-1} \Pempty{\uuu} (k! \Volume{\uuu})^{n-k+1}
            \One{k-\ell} (\uuu) \diff \uuu \diff r \diff z \diff L  \\
      &=  \tfrac{\density^{k+1}}{(k+1)!}
            \norm{\LGrass{k}{n}} \norm{\Region}
            \int\displaylimits_{r=0}^{r_0} r^{nk-1} e^{- \density r^n \nu_n}
            \!\!\! \int\displaylimits_{\uuu \in S^{k+1}} \!\!\!
              (k! \Volume{\uuu})^{n-k+1} \One{k-\ell}(\uuu) \diff \uuu \diff r,
    \label{eqn:terrible}
  \end{align}
  in which $S = S(L)$ is the $(k-1)$-dimensional unit sphere in $L$,
  as before.
  Making the substitution prepared in \eqref{eqn:gammasub}, we get
  \begin{align}
    \Expected{\ccon{\ell}{k}{}, r_0}
      &=  \density \norm{\LGrass{k}{n}} \norm{\Region}
          \tfrac{k!^{n-k+1}}{n (k+1)! \nu_n^k} \cdot \gamma (k, \density A_0)
          \int\displaylimits_{\uuu \in S^{k+1}}
              \Volume{\uuu}^{n-k+1} \One{k-\ell} (\uuu) \diff \uuu ,
  \end{align}
  and the claimed equation follows.
\eop

\ourparagraph{Constants.}
Note that Lemma \ref{lem:counting} implies Theorem \ref{thm:main}.
The integral in \eqref{eqn:counting} can be interpreted as a scaled
expectation of the volume of a $k$-simplex
spanned by $k+1$ points uniformly chosen on the unit sphere $\Sspace^{k-1}$.
With this in mind, we can split \eqref{eqn:counting} into a
monstrous but explicit constant factor and this expectation:
\begin{align}
  \Ccon{\ell}{k}{n} &=  \Factor{k}{n}
        \; \Expected{ \Volume{\uuu}^{n-k+1} \One{k-\ell} (\uuu) } ,
      \label{eqn:ClknExp}
\end{align}
which is over the $(k+1)$-multivariate uniform
distribution on the sphere, and the factor being
\begin{align}
  \Factor{k}{n} &= \tfrac{\sigma_n \cdot \sigma_{n-1} \cdot \ldots \cdot 
\sigma_{n-k+1}}
                 {\sigma_1 \cdot \sigma_2 \cdot \ldots \cdot \sigma_k}
           \tfrac{\Gamma(k) n^{k-1} k!^{n-k} \sigma_k^{k+1}}{(k+1) \sigma_n^k} ,
  \label{eqn:factor}
\end{align}
in which we use $\nu_n = \tfrac{\sigma_n}{n}$
to derive it from the form in \eqref{eqn:counting}.
To compute the coefficient for small values of $k$ and $n$,
it is helpful to recall that the measures of the unit spheres are
$\sigma_1 = 2$, $\sigma_2 = 2\pi$, $\sigma_3 = 4\pi$, $\sigma_4 = 2\pi^2$;
see Table \ref{tbl:g}.
\begin{table}[hbt]
  \centering
  \small{ \begin{tabular}{r||rrrr}
    \multicolumn{1}{c||}{$\Factor{k}{n}$}
      & $k=1$ & $2$ & $3$ & $4$   \\ \hline \hline
    $n = 2$  &  $1$  &  $\tfrac{4}{3} \pi$   &           &              \\
        $3$  &  $1$  &  $2 \pi^2$            &  $18 \pi$ &              \\
        $4$  &  $1$  &  $\tfrac{64}{3} \pi$  &  $1536$   &  $\tfrac{768}{5} 
\pi^2$
  \end{tabular} }
  \caption{Values of $\FactorOnly$ for small values of $k$ and $n$.}
  \label{tbl:g}
\end{table}
If we ignore the indicator function, this expectation is
the $(n-k+1)$-st moment of the $k$-dimensional volume of a random
inscribed $k$-simplex.
This allows us to write the number of $k$-simplices
that are upper bounds of intervals as
\begin{align}
  \sum_{\ell=0}^k \Ccon{\ell}{k}{n}
    &=  \Factor{k}{n} \; \Expected{ \Volume{\uuu}^{n-k+1} }.
  \label{eqn:UpperExp}
\end{align}
The value of this expectation was computed by Miles.
For $k = n$, this gives the number of top-dimensional Delaunay simplices:
\begin{align}
  \Dcon{n}{n}  &=  \sum_{\ell=0}^{n} \Ccon{\ell}{n}{n}
     =  \tfrac{2^{n+1} \pi^{(n-1)/2}}{n^2 (n+1)}
             \tfrac{\Gamma ( (n^2+1)/2 ) }
                   {\Gamma ( (n^2)/2   ) }
      \left[ \tfrac{\Gamma ( (n+2)/2   ) }
                   {\Gamma ( (n+1)/2   ) } \right]^n ;
\end{align}
see \cite{Mil71}.
Since every $n$-simplex has $n+1$ facets, and every $(n-1)$-simplex
in $\Rspace^n$
belongs to two $n$-simplices, we get $\Dcon{n-1}{n} = \tfrac{n+1}{2} 
\Dcon{n}{n}$.
This simple relation is perhaps more subtle than it appears.
It relies on Lemma \ref{lem:counting},
which writes the expected number of $(n-1)$-simplices in the Poisson--Delaunay
mosaic as a constant times $\density \norm{\Region}$
-- without lower order terms --
and on Lemma \ref{lem:boundarysize},
which implies that the difference between $K_0$ and $K_1$ is of lower order.
In dimension $n = 4$, we thus get
\begin{align}
  \Dcon{4}{4}  &=  \tfrac{32 \pi^{3/2}}{80}
             \tfrac{\Gamma (17/2) }
                   {\Gamma (8) }
      \left[ \tfrac{\Gamma (3) }
                   {\Gamma (5/2) } \right]^4 
          =  \tfrac{286}{9}
          =  31.77\ldots ,
    \label{eqn:S44}  \\
  \Dcon{3}{4}  &=  \tfrac{5}{2} \tfrac{286}{9}
          =  \tfrac{715}{9}
          =  79.44\ldots .
    \label{eqn:S34}
\end{align}

\section{Inscribed Simplices}
\label{sec:4}

In this section, we study the integral in \eqref{eqn:counting}
in more detail,
reducing it further to cones from the origin over the facets
of the simplex.

\ourparagraph{Volume decomposition.}
As before, we write $\uuu = ( u_0, u_1, \ldots, u_k )$ for a sequence
of $k+1$ distinct points in $\Sspace^{k-1}$ or, alternatively,
for the implied inscribed $k$-simplex.
For each $0 \leq i \leq k$, let $\uuu_i$ be the $k$-simplex
obtained by substituting $0$ for $u_i$,
and write $\Vcone{i}{} = \Volume{\uuu_i}$ for its $k$-dimensional volume.
Expressing the origin in terms of the points,
$0 = \sum_{i=0}^k \zeta_i u_i$ with $\sum_{i=0}^k \zeta_i = 1$,
we recall that the facet opposite to $u_i$ is visible from $0$
iff $\zeta_i < 0$.
Writing $\Signum{\zeta_i}$ for the sign of the $i$-th barycentric coordinate,
we therefore have
\begin{align}
  \Volume{\uuu}  &=  \sum_{i=0}^k \Signum{\zeta_i} \Vcone{i}{}.
  \label{eqn:VolumeSum}
\end{align}
The multiplicative group $\Zspace_2$ acts on the unit sphere by reflection.
This action is naturally extended to the action of $\Zspace_2^{k+1}$ on
$(k+1)$-tuples of points: for any vector $\ttt = (t_0, t_1, \ldots, t_k)$,
with $t_i \in \{-1, 1\}$ for $0 \leq i \leq k$,
we call $\ttt \uuu = ( t_0 u_0, t_1 u_1, \ldots, t_k u_k )$
the \emph{reflection} with \emph{signature} $\ttt$ of $\uuu$,
and we write $\Minus{\ttt}$ for the number of indices $i$ with $t_i = -1$.
Importantly, the reflection of a vertex does not affect the volume
of any cone.
We write $\Vsign{\ttt} = \Vsign{\ttt} (\uuu) = \sum_{i=0}^k t_i \Vcone{i}{}$
for the sum of positive and negative cone volumes.
Assuming $0$ is contained in the interior of the $k$-simplex $\uuu$,
the following lemma shows that it is the signed volume of $\ttt \uuu$.
\begin{lemma}[Volume Decomposition]
  \label{lem:voldec}
  Let $\uuu \in (\Sspace^{k-1})^{k+1}$ such that $0$ is contained
  in the interior of the $k$-simplex.
  Then $\Volume{\ttt \uuu}  =  \left| \Vsign{\ttt} (\uuu) \right|$,
  for every $\ttt \in \{-1, 1\}^{k+1}$.
\end{lemma}
\ourproof
  We reflect the vertices one by one to obtain $\ttt \uuu$ from $\uuu$
  and argue by induction on $\Minus{\ttt}$.
  By assumption, no facet of $\uuu$ is visible from $0$,
  so $\Volume{\uuu} = \sum_{i=0}^k \Vcone{i}{}$, which settles the base case.
  Assume without loss of generality that
  $\ttt = (-1, \ldots, -1, 1, \ldots, 1)$ with $\Minus{\ttt} = j$,
  and $\ttt^\prime = (-1, \ldots, -1, 1, \ldots, 1)$
  with $\Minus{\ttt^\prime} = j-1$.
  By induction, the volume of $\ttt^\prime \uuu$
  is $\pm \Vsign{\ttt^\prime} (\uuu)$, i.e., either
  $\Volume{\ttt^\prime \uuu} = 
-\Vcone{0}{}-\ldots-\Vcone{j-1}{}+\Vcone{j}{}+\ldots+\Vcone{k}{}$ or
  $\Volume{\ttt^\prime \uuu} =  
\Vcone{0}{}+\ldots+\Vcone{j-1}{}-\Vcone{j}{}-\ldots-\Vcone{k}{}$,
  depending on which of the two expressions is positive.
  Reflecting $u_j$ either changes the orientation of the inscribed $k$-simplex,
  or it does not.
  In case the orientation is changed,
  the reflection changes the visibility of exactly one facet,
  namely the one opposite to $u_j$,
  and by \eqref{eqn:VolumeSum} we get either
  $\Volume{\ttt \uuu} = 
-\Vcone{0}{}-\ldots-\Vcone{j-1}{}-\Vcone{j}{}+\Vcone{j+1}{}+\ldots+\Vcone{k}{}$ 
or
  $\Volume{\ttt \uuu} =  
\Vcone{0}{}+\ldots+\Vcone{j-1}{}+\Vcone{j}{}-\Vcone{j+1}{}-\ldots-\Vcone{k}{}$.
  In case the orientation is preserved, the reflection changes the visibility
  of every facet but one, namely the one opposite to $u_j$,
  and again by \eqref{eqn:VolumeSum} we get either
  $\Volume{\ttt \uuu} =  
\Vcone{0}{}+\ldots+\Vcone{j-1}{}+\Vcone{j}{}-\Vcone{j+1}{}-\ldots-\Vcone{k}{}$ 
or
  $\Volume{\ttt \uuu} = 
-\Vcone{0}{}-\ldots-\Vcone{j-1}{}-\Vcone{j}{}+\Vcone{j+1}{}+\ldots+\Vcone{k}{}$.
  Summing up, in all cases we have $\Volume{\ttt \uuu} = |\Vsign{\ttt} (\uuu)|$.
\eop

\ourparagraph{Visibility.}
There are several useful consequences of Lemma \ref{lem:voldec},
which we now state.
For example, for almost every inscribed $k$-simplex, there are
precisely two signatures for which the corresponding
reflections produce a $k$-simplex that contain the origin.
To produce one, we reflect every vertex opposite a facet visible
from $0$,
and to produce the other, we reflect every vertex opposite a
facet that is not visible from $0$.
If the first simplex corresponds to $\ttt = (t_0, t_1, \ldots, t_k)$,
then the second corresponds to $- \ttt = (-t_0, -t_1, \ldots, -t_k)$,
which we refer to as the \emph{complementary} signature.
\begin{corollary}[Reflections and Visibility]
  \label{cor:refvis}
  Let $\uuu \in (\Sspace^{k-1})^{k+1}$ such that $0$ is contained
  in the interior of the $k$-simplex, and let $\ttt \in \{-1,1\}^{k+1}$.
  \begin{enumerate}
    \item[1.] After reflecting a subset of the vertices, the visible facets
      are either the ones opposite to the reflected vertices, or all others.
      Specifically, if $\Vsign{\ttt} (\uuu) > 0$, then there are $\Minus{\ttt}$
      visible facets, each one opposite a reflected vertex,
      and if $\Vsign{\ttt} (\uuu) < 0$, then there are $k - \Minus{\ttt} + 1$
      visible facets, each one opposite a non-reflected vertex.
    \item[2.] The simplices $\ttt \uuu$ and $- \ttt \uuu$
      are central reflections of each other; in particular,
      they have the same volume and corresponding facets are visible from $0$.
  \end{enumerate}
\end{corollary}

Fact 1 in Corollary \ref{cor:refvis} is a direct consequence
of \eqref{eqn:VolumeSum} and Lemma \ref{lem:voldec},
and Fact 2 is clear for geometric reasons.
The following simple facts will be useful in our computations.
\begin{lemma}[Visibility of Facets]
  \label{lem:visfacets}
  Let $\uuu \in (\Sspace^{k-1})^{k+1}$ such that $0$ is contained
  in the interior of the $k$-simplex, and let $\ttt \in \{-1,1\}^{k+1}$.
  \begin{enumerate}
    \item[1.] The origin, $0$, is contained in the interior of the $k$-simplex
      $\ttt \uuu$ iff $\Minus{\ttt} = 0$ or $k+1$.
    \item[2.] $\Minus{\ttt} = 0$ implies $\Vsign{\ttt} (\uuu) > 0$ and, 
equivalently,
      $\Minus{\ttt} = k+1$ implies $\Vsign{\ttt} (\uuu) < 0$.
    \item[3.] If a set of facets is visible from $0$ for $\ttt \uuu$,
      then there is no signature $\ttt'$ such that the complementary set of 
facets
      is visible from $0$ for $\ttt' \uuu$.
    \item[4.] $\Minus{\ttt} = 1$ implies $\Vsign{\ttt} (\uuu) > 0$ and, 
equivalently,
              $\Minus{\ttt} = k$ implies $\Vsign{\ttt} (\uuu) < 0$.
  \end{enumerate}
\end{lemma}
\ourproof
  By assumption on $\uuu$, the only signatures for which all terms $t_i 
\Vcone{i}{}$ have
  the same sign are the ones for which $\Minus{\ttt} = 0$ or $\Minus{\ttt} = 
k+1$.
  Fact 1 follows and implies Fact 2.

  To see Fact 3, we express $\Volume{\ttt \uuu}$ using \eqref{eqn:VolumeSum},
  getting a negative coefficient for every visible facet.
  Nevertheless, the sum of signed cone volumes is positive.
  If the visibility of all facets could be reversed, \eqref{eqn:VolumeSum}
  would give a negative volume, which is a contradiction.
  Fact 4 follows:  it is possible to see every single facet of the simplex,
  while seeing none of the others,
  hence it is impossible to see the complementary $k$ facets together from $0$.
\eop

Fact 1 of Lemma \ref{lem:visfacets} appears already in Wendel \cite{Wen62},
who generalized it to compute the probability that all points of a finite
set sampled independently and uniformly on a sphere lie inside a hemisphere.

\ourparagraph{Spherical expectations.}
We now return to \eqref{eqn:ClknExp}.
The probability space for it is
$\mathcal{P} = (\{\Sspace^{k-1}, \mathfrak{U}\})^{\otimes(k+1)}$,
in which $\mathfrak{U}$ is the uniform measure on the sphere,
and the random variables $u_i$ are just the projections. 
Note that every inscribed simplex, $\uuu$, corresponds to a unique point
configuration $\bar{\uuu} \in \RPspace^{k-1}$ obtained by projecting
$u_i$ from the sphere to the projective space.
Likewise, every $k$-simplex with vertices in $\RPspace^{k-1}$
corresponds to $2^{k+1}$ $k$-simplices inscribed in $\Sspace^{k-1}$.
It allows us to decompose the probability space as
$\mathcal{P} = (\{\RPspace^{k-1}, \mathfrak{U}^\prime\} \otimes \{\{-1,1\},
 \mathfrak{B}\})^{\otimes(k+1)}$,
in which $\mathfrak{U}^\prime$ is the uniform measure on the projective space
and $\mathfrak{B}$ is the uniform measure on $\Zspace_2$.
In other words, we decompose the uniform measure on the sphere
as the measure on orbits under the action of $\Zspace_2^{k+1}$
times the Haar measure on the group.
Write $\EE_\uuu$ for the expectation taken over the sphere,
$\EE_{\bar{\uuu}}$ for the expectation over the projective space
and $\EE_{\bar{\uuu},\ttt}$ for the expectation over the projective space and 
the group.
We use the probabilistic formalism only locally,
to decompose the expectation in \eqref{eqn:ClknExp} further into
expectations involving volumes of cones.
We recall that the volume of $\ttt \uuu$ is either $\Vsign{\ttt} (\uuu)$
or $- \Vsign{\ttt} (\uuu)$.
For each $0 \leq \ell \leq k \leq n$, we write the expectation in 
\eqref{eqn:ClknExp} as
\begingroup
\allowdisplaybreaks
\begin{align}
  \EAux{\ell}{k}{n}  &= \EE_\uuu [ \Volume{\uuu} ^{n-k+1} \One{k-\ell} (\uuu) ]
    \label{eqn:reduction0} \\
                     &= \EE_{\bar{\uuu},\ttt}
                        [| \Vsign{\ttt} (\bar{\uuu}) |^{n-k+1} \One{k-\ell} 
(\bar{\uuu}, \ttt) ]
    \label{eqn:reduction1} \\
                     &=  \tfrac{1}{2^{k+1}} \!\! \sum_{\Minus{\ttt} = k-\ell}
                         \EE_{\bar{\uuu}}{[|\Vsign{\ttt} (\bar{\uuu})|^{n-k+1} 
\One{\Vsign{\ttt} (\bar{\uuu}) >0}]}
                      +  \tfrac{1}{2^{k+1}} \!\! \sum_{\Minus{\ttt} = \ell+1}
                         \EE_{\bar{\uuu}}{[|\Vsign{\ttt} (\bar{\uuu})|^{n-k+1} 
\One{\Vsign{\ttt} (\bar{\uuu}) <0}]}
  \label{eqn:reduction2} \\
                     &=  \tfrac{1}{2^k}
                         \sum_{\Minus{\ttt} = k-\ell}
                         \EE_{\bar{\uuu}}{[\Vsign{\ttt} (\bar{\uuu})^{n-k+1} 
\One{\Vsign{\ttt} (\bar{\uuu}) >0}]} 
  \label{eqn:reduction3} \\
                     &=  \tfrac{1}{2^k} \tbinom{k+1}{k-\ell}
                         \; \EE_{\uuu}{[\Vsign{\ttt_{k-\ell}} (\uuu)^{n-k+1} 
                         \One{\Vsign{\ttt_{k-\ell}} (\uuu) > 0}]} ,
  \label{eqn:reduction4}
\end{align}
\endgroup
in which $\ttt_{k-\ell}$ in \eqref{eqn:reduction4} is an arbitrary signature
with $\Minus{\ttt} = k-\ell$.
The transition to \eqref{eqn:reduction1} is possible because
for a fixed $\ttt$, $\Vsign{\ttt}$ is the same for all simplices in an orbit,
and the transition to \eqref{eqn:reduction2} is justified
by the first fact in Corollary \ref{cor:refvis}.
We get \eqref{eqn:reduction3} by observing that the two sums in
\eqref{eqn:reduction2} are over complementary signatures,
and we get \eqref{eqn:reduction4} because relabeling the vertices
does not change the expected volume.
We can remove the bar in the last transition again because $\Vsign{\ttt_{k-\ell}}$
is the same along the orbits.
Combining \eqref{eqn:ClknExp} with \eqref{eqn:reduction4},
we finally get the expected number of intervals of each type in terms
of the spherical expectations.
\begin{lemma}[Reduction to Spherical Expectations]
  \label{lem:reduction}
  For every $0 \leq \ell \leq k \leq n$, we have
  $\Ccon{\ell}{k}{n} = \Factor{k}{n} \cdot \EAux{\ell}{k}{n}$.
\end{lemma}

\ourparagraph{Moments.}
We present results on random simplices and their volumes
needed to derive the desired expectations.
The first result gives the moments of the cone volumes;
see \cite[Equation 2.11]{Mil69} and with a minor correction
\cite[Equation (23)]{Mil71}.
Let $u_1, u_2, \ldots, u_k$ be independently and uniformly
distributed on $\Sspace^{n-1}$,
and write $\uuu_0 = (0, u_1, \ldots, u_k)$ for the corresponding cone,
which is a $k$-simplex.
Then for any integer $a \geq 0$, the expectation of $\Volume{\uuu_0}^a$ is
\begin{align}
  \Mnt{k}{n}{a}  &=  \frac{1}{k!^a}
    \left[ \tfrac{\Gamma (  n   /2 ) }
                 {\Gamma ( (n+a)/2 ) } \right]^{k-1}
    \prod_{i=1}^{k-1}
           \tfrac{\Gamma ( (n-k+a+i)/2 )}
                 {\Gamma ( (n-k+i)/2   )} .
  \label{eqn:moments}
\end{align}
Besides these moments, we also need the mixed
moments to get our results.
At the time of writing, we have these only for pairs of cones.
Given $\uuu = (u_0, u_1, \ldots, u_n) \in (\Sspace^{n-1})^{n+1}$,
we recall that $\uuu_i$ is the cone obtained by substituting $0$
for $u_i$.
\begin{lemma}[Pairwise Mixed Moments]
  \label{lem:mixedmoments}
  Let $\uuu$ be a sequence of $n+1$ independently and uniformly
  distributed points on $\Sspace^{n-1}$.
  Then for any $0 \leq i < j \leq n$ and integers $a, b \geq 0$,
  the expectation of $\Vcone{i}{a} \Vcone{j}{b}$ is
  \begin{align}
    \xMnt{n}{a}{b}
       &=  \tfrac{ \Mnt{n-1}{n}{a+b} }{ n^{a+b} }
           \left[ \tfrac{ \Gamma(n/2)     }{ \Gamma(1/2)     } \right]^2
                  \tfrac{ \Gamma((a+1)/2) }{ \Gamma((n+a)/2) }
                  \tfrac{ \Gamma((b+1)/2) }{ \Gamma((n+b)/2) } .
    \label{eqn:mixedmoments}
  \end{align}
\end{lemma}

\ourproof
  Note that $\Vcone{i}{} = \tfrac{1}{n} h_i A$
  and $\Vcone{j}{} = \tfrac{1}{n} h_j A$,
  in which $A$ is the $(n-1)$-dimensional volume of the shared facet
  of $\uuu_i$ and $\uuu_j$, and $h_i, h_j$ are the distances
  of the points $u_i, u_j$ from the hyperplane spanned by the shared facet.
  For geometric reasons, it is clear that $h_i, h_j, A$
  are independent; see \cite{Mil69} for details.
  Hence, we get
  \begin{align}
    \Expected{\Vcone{i}{a} \Vcone{j}{b}} &= \tfrac{1}{n^{a+b}}
     \; \Expected{h_i^a} \; \Expected{h_j^b} \; \Expected{A^{a+b}} ,
  \end{align}
  with $\Expected{A^{a+b}} = \Mnt{n-1}{n}{a+b}$ by \eqref{eqn:moments}.
  The value for $\Expected{h_i^a}$ given in \cite{Mil69},
  right before Formula (2.11), is
  $\Gamma(\tfrac{n}{2}) / \Gamma(\tfrac{1}{2})$ times
  $\Gamma(\tfrac{a+1}{2}) / \Gamma(\tfrac{n+a}{2})$.
  Substituting the analogous expression for $\Expected{h_j^b}$
  gives the claimed relation.
\eop

We illustrate \eqref{eqn:moments} and \eqref{eqn:mixedmoments}
by computing $\MntOnly$ and $\xMntOnly$ for a selected set of
small values of $k, n, a, b$, chosen so the results will be
useful in Section \ref{sec:5}.
\begin{table}[hbt]
  \centering
  \small{ \begin{tabular}{r||ccc}
                                            &       &         &         \\
    \multicolumn{1}{c||}{$\Mnt{k}{n}{a}$} & $a=1$ & $a = 2$ & $a = 3$ \\ \hline 
\hline
    $k = n = 2$ & $\tfrac{1}{\pi}$      & $\tfrac{1}{8}$   & $\tfrac{1}{6 \pi}$ 
\\
            $3$ & $\tfrac{\pi}{48}$     & $\tfrac{1}{162}$ &                    
\\
            $4$ & $\tfrac{8}{81 \pi^2}$ &                  &
  \end{tabular} }
  \hspace{0.15in}
  \small{ \begin{tabular}{r||ccc}
                                            & $a=1$ & $a=2$ & $a=2$ \\
    \multicolumn{1}{c||}{$\xMnt{n}{a}{b}$} & $b=1$ & $b=1$ & $b=2$ \\ \hline 
\hline
    $n = 2$ & $\tfrac{1}{\pi^2}$  & $\tfrac{1}{8 \pi}$   &                    \\
        $3$ & $\tfrac{1}{216}$    &                      &                    \\
        $4$ &                     &                      &
  \end{tabular} }
  \caption{Values of $\MntOnly$ for small values of $k, n, a$ on the
    \emph{left}, and values of $\xMntOnly$ for small values of $n, a, b$
    on the \emph{right}.}
  \label{tbl:mumm}
\end{table}

\section{Constants}
\label{sec:5}

Being done with the general facts, we apply them to give concrete expressions
for the expected numbers of intervals of the radius function in dimensions up 
to 
four.
We will mainly compute expectations by integrations --
using Equation \eqref{eqn:counting} and Lemma \ref{lem:reduction} --
and rarely resort to the linear relations that connect the expectations.

\subsection{Two Dimensions}
\label{sec:51}

As a warm-up exercise, we begin with a Poisson point process in $\Rspace^2$.
We have $\Ccon{0}{0}{2} = 1$ and $\Ccon{0}{k}{2} = 0$ for $k > 0$ because
all vertices are critical.
To compute the remaining constants, we need the spherical expectations
given in \eqref{eqn:reduction4}:
\begin{align}
  2            \EAux{1}{1}{2}  &=  \Expected{(\Vcone{0}{} + \Vcone{1}{})^2}
                                =  2 \Expected{\Vcone{0}{2}} + 2 
\Expected{\Vcone{0}{} \Vcone{1}{}}, 
    \label{eqn:firstExp2} 
\end{align}
in which we get the right-hand side because expectations do not
change under re-indexing.
The expectation is with respect to the uniform distribution on $\Sspace^0$,
which is a pair of points.
We have $\Vcone{0}{} = \Vcone{1}{} = 1$ and therefore $2 \EAux{1}{1}{2} = 4$.
We also need
\begin{align}
  \tfrac{4}{3} \EAux{1}{2}{2}  &=  \Expected{\Vcone{0}{} + \Vcone{1}{} - 
\Vcone{2}{}}
                                =  \Expected{\Vcone{0}{}} ,
    \label{eqn:secondExp2} \\
  4            \EAux{2}{2}{2}  &=  \Expected{\Vcone{0}{} + \Vcone{1}{} + 
\Vcone{2}{}}
                                =  3 \Expected{\Vcone{0}{}} ,
    \label{eqn:thirdExp2}
\end{align}
which both satisfy $\Minus{\ttt} \leq 1$, as did \eqref{eqn:firstExp2},
so Lemma \ref{lem:visfacets} applies and we can remove the indications,
which we did.
These two expectations are with respect to the uniform distribution
on $\Sspace^1$.
Using \eqref{eqn:moments} to compute $\Expected{\Vcone{0}{}}$, we get
$\tfrac{4}{3} \EAux{1}{2}{2} = \Mnt{2}{2}{1} = \tfrac{1}{\pi}$,
and similarly $4 \EAux{2}{2}{2} = 3 \Mnt{2}{2}{1} = \tfrac{3}{\pi}$.
Retrieving $\Factor{1}{2} = 1$ and $\Factor{2}{2} = \tfrac{4 \pi}{3}$
from Table \ref{tbl:mumm},
we can now use Lemma \ref{lem:reduction} to get the corresponding constants:
\begin{align}
  \Ccon{1}{1}{2}  &=  \Factor{1}{2} \cdot \EAux{1}{1}{2}
                   =  1 \cdot \tfrac{1}{2} \cdot 4
                   =  2 ,
    \label{eqn:C112} \\
  \Ccon{1}{2}{2}  &=  \Factor{2}{2} \cdot \EAux{1}{2}{2}
                   =  \tfrac{4 \pi}{2} \cdot \tfrac{3}{4} \cdot \tfrac{1}{\pi}
                   =  1 ,
    \label{eqn:C122} \\
  \Ccon{2}{2}{2}  &=  \Factor{2}{2} \cdot \EAux{2}{2}{2}
                   =  \tfrac{4 \pi}{2} \cdot \tfrac{1}{4} \cdot 3 \tfrac{1}{\pi}
                   =  1 .
    \label{eqn:C222}
\end{align}
This justifies the entries of the left matrix in Table \ref{tbl:constants}.
Note that $\Ccon{0}{0}{2} - \Ccon{1}{1}{2} + \Ccon{2}{2}{2} = 0$,
which agrees with the discrete Morse relation stated as Lemma 
\ref{lem:discreteMorse}.
Indeed, it makes sense to use this relation as a check of correctness as we
have refrained from using it during the derivation of the constants.

\vspace{0.1in} \noindent {\sc Remark.}
As pointed out by G\"{u}nter Rote, the computations
for the critical edges generalize to $n$ dimensions.
Indeed, $\Factor{1}{n} = 1$ and
$2 \EAux{1}{1}{2} = \Expected{(\Vcone{0}{} + \Vcone{1}{})^n} = 2^n$, which gives
\begin{align}
  \Ccon{1}{1}{n}  &=  \Factor{1}{n} \cdot \EAux{1}{1}{n}
                   =  2^{n-1}.
  \label{rem:edges}
\end{align}

%
\ourparagraph{Simplices in the Poisson--Delaunay mosaic.}
%
For completeness, we also compute the expected numbers of simplices in
the $2$-dimensional Poisson--Delaunay mosaic, which are of course known:
\begin{align}
  \Dcon{0}{2}  &=  \Ccon{0}{0}{2}  =  1 ,
    \label{eqn:D02} \\
  \Dcon{1}{2}  &=  \Ccon{1}{1}{2} + \Ccon{1}{2}{2} = 3 ,
    \label{eqn:D12} \\
  \Dcon{2}{2}  &=  \Ccon{1}{2}{2} + \Ccon{2}{2}{2} = 2 .
    \label{eqn:D22}
\end{align}
We have $\Dcon{0}{2} - \Dcon{1}{2} + \Dcon{2}{2} = 0$,
which is consistent with the Euler relation in the plane.
Note that \eqref{eqn:C222} and \eqref{eqn:D22} imply that about half the
Delaunay triangles are critical.
The geometric reason behind this fact is an observation by Miles \cite{Mil70}
that a Delaunay triangle is acute with probability $\tfrac{1}{2}$.

\subsection{Three Dimensions}
\label{sec:52}

We have $\Ccon{0}{0}{3} = 1$ and $\Ccon{0}{k}{3} = 0$ for $k > 0$
because every vertex is critical,
and we know $\Ccon{1}{1}{3} = 4$ for the critical edges from \eqref{rem:edges}.
To compute the remaining constants in $\Rspace^3$,
we need some spherical expectations:
\begin{align}
  \tfrac{4}{3} \EAux{1}{2}{3}  &=  \Expected{(\Vcone{0}{} + \Vcone{1}{} - 
\Vcone{2}{})^2} 
                                =  3 \Expected{\Vcone{0}{2}} - 2 
\Expected{\Vcone{0}{} \Vcone{1}{}},    \\
  4            \EAux{2}{2}{3}  &=  \Expected{(\Vcone{0}{} + \Vcone{1}{} + 
\Vcone{2}{})^2}            
                                =  3 \Expected{\Vcone{0}{2}} + 6 
\Expected{\Vcone{0}{} \Vcone{1}{}},    
\end{align}
in which the expectations are with respect to the uniform distribution on the 
circle.
We get $\Expected{\Vcone{0}{2}} = \Mnt{2}{2}{2} = \tfrac{1}{8}$
from \eqref{eqn:moments} and
$\Expected{\Vcone{0}{} \Vcone{1}{}} = \xMnt{2}{1}{1} = \tfrac{1}{\pi^2}$
from \eqref{eqn:mixedmoments}; see also Table \ref{tbl:mumm}.
Using again Lemma \ref{lem:visfacets} to omit indicators, we furthermore have
\begin{align}
  2            \EAux{2}{3}{3}  &=  \Expected{\Vcone{0}{} + \Vcone{1}{} + 
\Vcone{2}{} - \Vcone{3}{}}
                                =  2 \Expected{\Vcone{0}{}},                     
        \\
  8            \EAux{3}{3}{3}  &=  \Expected{\Vcone{0}{} + \Vcone{1}{} + 
\Vcone{2}{} + \Vcone{3}{}}
                                =  4 \Expected{\Vcone{0}{}} ,
\end{align}
in which the expectations are with respect to the uniform distribution
on the $2$-dimensional sphere.
For the moment we skip the computation
of $\tfrac{8}{6} \EAux{1}{3}{3} = \Expected{|\Vcone{0}{} + \Vcone{1}{} - 
\Vcone{2}{} - \Vcone{3}{}|}$.
We get $\Expected{\Vcone{0}{}} = \Mnt{3}{3}{1} = \tfrac{\pi}{48}$ from 
\eqref{eqn:moments}.
Multiplying the spherical expectation with the corresponding factors
in Lemma \ref{lem:reduction}, we get the corresponding entries of the
middle matrix in Table \ref{tbl:constants}:
\begin{align}
  \Ccon{1}{2}{3}  &=  \Factor{2}{3} \cdot \EAux{1}{2}{3}
             =  2 \pi^2 \cdot \tfrac{3}{4} \cdot (3 \tfrac{1}{8} - 2 
\tfrac{1}{\pi^2})
             =  \tfrac{9}{16}\pi^2 - 3 
             =  2.55\ldots ,                                                     
      \\
  \Ccon{2}{2}{3}  &=  \Factor{2}{3} \cdot \EAux{2}{2}{3}
             =  2 \pi^2 \cdot \tfrac{1}{4} \cdot (3 \tfrac{1}{8} + 6 
\tfrac{1}{\pi^2})
             =  \tfrac{3}{16}\pi^2 + 3 
             =  4.85\ldots ,                                                      
      \\
  \Ccon{2}{3}{3}  &=  \Factor{3}{3} \cdot \EAux{2}{3}{3}
             =  18 \pi  \cdot \tfrac{1}{2} \cdot 2 \tfrac{\pi}{48}
             =  \tfrac{3}{8}\pi^2
             =  3.70\ldots ,                                                    
      \\
  \Ccon{3}{3}{3}  &=  \Factor{3}{3} \cdot \EAux{3}{3}{3}
             =  18 \pi  \cdot \tfrac{1}{8} \cdot 4 \tfrac{\pi}{48}
             =  \tfrac{3}{16}\pi^2
             =  1.85\ldots .
\end{align}
We can compute the remaining $\Ccon{1}{3}{3}$ either by Euler formula
or from \eqref{eqn:UpperExp}, which gives the constant in the number of 
$3$-simplices
in the Poisson--Delaunay mosaic as $\Dcon{3}{3} = \tfrac{24}{35}\pi^2$;
see also \cite{ScWe08}.
This gives
\begin{align}
  \Ccon{1}{3}{3}  &=  \tfrac{69}{560}\pi^2 
             =  1.21\ldots ,
\end{align}
which completes the justification of the entries of the middle matrix
in Table \ref{tbl:constants}.
We use Lemma \ref{lem:discreteMorse} to check the numbers of
critical simplices and get
$\Ccon{0}{0}{3} - \Ccon{1}{1}{3} + \Ccon{2}{2}{3} - \Ccon{3}{3}{3} = 0$, as 
required.

\ourparagraph{Simplices in the Poisson--Delaunay mosaic.}
While the expected numbers of simplices in the Poisson--Delaunay mosaic
in $\Rspace^3$ are known \cite{ScWe08}, it is easy to compute them from
the above constants:
\begin{align}
  \Dcon{0}{3}  &=  \Ccon{0}{0}{3}
          =  1 ,                                                  \\
  \Dcon{1}{3}  &=  \Ccon{1}{1}{3} + \Ccon{1}{2}{3} + \Ccon{1}{3}{3}
          =  \tfrac{24}{35} \pi^2 + 1
          =  7.76\ldots,                                          \\
  \Dcon{2}{3}  &=  \Ccon{1}{2}{3} + \Ccon{2}{2}{3} + 2\Ccon{1}{3}{3} + 
\Ccon{2}{3}{3}
          =  \tfrac{48}{35} \pi^2
          =  13.53\ldots,                                         \\
  \Dcon{3}{3}  &=  \Ccon{1}{3}{3} + \Ccon{2}{3}{3} + \Ccon{3}{3}{3}
          =  \tfrac{24}{35} \pi^2
          =  6.76\ldots  .
\end{align}
This completes the entries in the second row of Table \ref{tbl:PD-constants}.
As a final check of correctness, we compute the alternating sum,
which gives $\Dcon{0}{3} - \Dcon{1}{3} + \Dcon{2}{3} - \Dcon{3}{3} = 0$, as 
required.

\subsection{Four Dimensions}
\label{sec:53}

In four dimensions, we compute most of the constants directly,
but use knowledge of $\Dcon{4}{4}$ and $\Dcon{3}{4}$ to get $\Ccon{1}{4}{4}$
and $\Ccon{2}{4}{4}$.
We have $\Ccon{0}{0}{4} = 1$ and $\Ccon{0}{k}{4} = 0$ for $k > 0$
because every vertex is critical,
and $\Ccon{1}{1}{4} = 8$ by \eqref{rem:edges},
so we proceed to the remaining constants.

\ourparagraph{Triangles as upper bounds.}
Here we count the critical triangles and edge-triangle pairs.
Starting with $\Ccon{1}{2}{4}$, we have $\Minus{\ttt} = 1$ reflection,
and by Lemma \ref{lem:visfacets} this implies $\Vsign{\ttt} > 0$.
We therefore get
\begin{align}
  \tfrac{4}{3} \EAux{1}{2}{4}
    &=    \Expected{(\Vcone{0}{} + \Vcone{1}{} - \Vcone{2}{})^3}  \\
    &=    \Expected{\Vcone{0}{3} \!+\! \Vcone{1}{3} \!-\! \Vcone{2}{3}
                + 3 (\Vcone{0}{2} \Vcone{1}{}
              \!-\!  \Vcone{0}{2} \Vcone{2}{}
              \!+\!  \Vcone{1}{2} \Vcone{0}{}
              \!-\!  \Vcone{1}{2} \Vcone{2}{}
              \!+\!  \Vcone{2}{2} \Vcone{0}{}
              \!+\!  \Vcone{2}{2} \Vcone{1}{})
                - 6  \Vcone{0}{}  \Vcone{1}{}  \Vcone{2}{} } \\
    &=    \Expected{\Vcone{0}{3}}
      + 6 \Expected{\Vcone{0}{2} \Vcone{1}{}}
      - 6 \Expected{\Vcone{0}{} \Vcone{1}{} \Vcone{2}{}} .
   \label{eqn:E124-3}
\end{align}
From \eqref{eqn:moments} and \eqref{eqn:mixedmoments} we get
$\Expected{\Vcone{0}{3}} = \Mnt{2}{2}{3} = \tfrac{1}{6\pi}$ and
$\Expected{\Vcone{0}{2} \Vcone{1}{}} = \xMnt{2}{2}{1} = \tfrac{1}{8\pi}$.
Note that $\Vcone{0}{}$ and $\Vcone{1}{}$ are independent
in two dimensions, so we also have
$\Expected{\Vcone{0}{2} \Vcone{1}{}} =
  \Expected{\Vcone{0}{2}} \; \Expected{\Vcone{1}{}} =
  \Mnt{2}{2}{2} \; \Mnt{2}{2}{1}$,
which gives the same result.
For the remaining term, we need a convenient description of the three points
uniformly chosen on the unit circle.
Fixing $u_0$, we parametrize $u_1$ and $u_2$ by the angles
$\alpha, \beta \in [-\pi, \pi]$ they form with $u_0$.
In this setup, we have
$\Vcone{0}{} = \tfrac{1}{2} | \sin (\alpha - \beta) |$,
$\Vcone{1}{} = \tfrac{1}{2} | \sin \beta |$,
$\Vcone{2}{} = \tfrac{1}{2} | \sin \alpha |$,
where $\alpha$ and $\beta$ are uniformly distributed over $[-\pi, \pi]$. We 
notice that this also implies that $\Vcone{i}{}$ and $\Vcone{j}{}$ are 
independent whenever $i \neq j$.
The moment can now be computed as
\begin{align}
  \Expected{\Vcone{0}{} \Vcone{1}{} \Vcone{2}{}}
    &= \tfrac{1}{8}\;\Expected{|\sin \alpha| |\sin \beta| |\sin 
(\alpha-\beta)|} 
  \label{eqn:tripleMoment1} \\
    &= \tfrac{1}{8} \tfrac{1}{4 \pi^2}
       \int\displaylimits_{\alpha=-\pi}^\pi
       \int\displaylimits_{\beta =-\pi}^\pi
        |\sin \alpha| |\sin \beta| |\sin (\alpha-\beta)| \diff \alpha \diff 
\beta 
  \label{eqn:tripleMoment2} \\
    &= \tfrac{1}{8 \pi^2}
       \int\displaylimits_{\alpha=0}^\pi
       \int\displaylimits_{\beta =0}^\pi
        \sin \alpha \sin \beta |\sin (\alpha-\beta)| \diff \alpha \diff \beta ,
  \label{eqn:tripleMoment3}
\end{align}
in which \eqref{eqn:tripleMoment3} is true because the expression does not 
change
under transformations $\alpha \mapsto \alpha + \pi$ and $\beta \mapsto \beta + 
\pi$.
Computing the integral either by splitting cases or using any mathematical
software, we see that the moment evaluates to $\tfrac{3}{32 \pi}$.
Next, we proceed to the critical triangles, computing $\Ccon{2}{2}{4}$.
For this, we need
\begin{align}
  4 E_{22}^4  &=  \Expected{ (\Vcone{0}{} + \Vcone{1}{} + \Vcone{2}{})^3 } 
               =  3 \Expected{\Vcone{0}{3}}
               + 18 \Expected{\Vcone{0}{2} \Vcone{1}{}}
               +  6 \Expected{\Vcone{0}{} \Vcone{1}{} \Vcone{2}{}} .
\end{align}
Plugging these results into Lemma \ref{lem:reduction}, we get
\begin{align}
  \label{eqn:ans43}
  \Ccon{1}{2}{4}  &=  \Factor{2}{4} \cdot E_{12}^4
             =  \tfrac{64 \pi}{3} \cdot \tfrac{3}{4} \cdot
                ( \tfrac{1}{6\pi} + 6 \tfrac{1}{8\pi} - 6 \tfrac{3}{32\pi} ) 
             =  \tfrac{17}{3}
             =  5.66\ldots , \\
  \Ccon{2}{2}{4}  &=  \Factor{2}{4} \cdot E_{22}^4
             =  \tfrac{64 \pi}{3} \cdot \tfrac{1}{4} \cdot
                ( 3 \tfrac{1}{6\pi} + 18 \tfrac{1}{8\pi} + 6 \tfrac{3}{32\pi} )
             =  \tfrac{53}{3}
             =  17.66\ldots .
\end{align}

\ourparagraph{Tetrahedra as upper bounds.}
Here we count the critical tetrahedra,
triangle-tetrahedron pairs, and edge-tetrahedron quadruplets.
Starting with $\Ccon{1}{3}{4}$, we need the second moment of the volumes
of cones with two visible facets.
Setting $\ttt = (1, 1, -1, -1)$ and recalling that $- \ttt = (-1, -1, 1, 1)$, 
we 
get
\begin{align}
  \tfrac{4}{3} E_{13}^4
    &=  \Expected{ (\Vcone{0}{} + \Vcone{1}{} - \Vcone{2}{} - \Vcone{3}{})^2
                   \One{\Vsign{\ttt} > 0} } \\
    &=  \tfrac{1}{2} \left(
        \Expected{ (\Vcone{0}{} + \Vcone{1}{} - \Vcone{2}{} - \Vcone{3}{})^2
                   \One{\Vsign{\ttt} > 0} } 
      + \Expected{ (- \Vcone{0}{} - \Vcone{1}{} + \Vcone{2}{} + \Vcone{3}{})^2
                   \One{\Vsign{-\ttt} > 0} } \right) \\
    &=  \tfrac{1}{2} \left(
        \Expected{ (\Vcone{0}{} + \Vcone{1}{} - \Vcone{2}{} - \Vcone{3}{})^2
                   \One{\Vsign{\ttt} > 0} } 
      + \Expected{ (\Vcone{0}{} + \Vcone{1}{} - \Vcone{2}{} - \Vcone{3}{})^2
                   \One{\Vsign{\ttt} < 0} } \right) \\
    &=  \tfrac{1}{2} \left(
        \Expected{ (\Vcone{0}{} + \Vcone{1}{} - \Vcone{2}{} - \Vcone{3}{})^2 }
                     \right) \\
    &=  2 \Expected{\Vcone{0}{2}} - 2 \Expected{\Vcone{0}{} \Vcone{1}{}}.
\end{align}
We get $\Expected{\Vcone{0}{2}} = \Mnt{3}{3}{2} = \tfrac{1}{162}$
from \eqref{eqn:moments}, and
$\Expected{\Vcone{0}{} \Vcone{1}{}} = \xMnt{3}{1}{1} = \tfrac{1}{216}$
from \eqref{eqn:mixedmoments}.
Moving on to $\Ccon{2}{3}{4}$ and to $\Ccon{3}{3}{4}$, we need
\begin{align}
  2 E_{23}^4
    &=  \Expected{ (\Vcone{0}{} + \Vcone{1}{} + \Vcone{2}{} - \Vcone{3}{})^2 }
     =  4 \Expected{ \Vcone{0}{2} },                                         \\
  8 E_{33}^4
    &=  \Expected{ (\Vcone{0}{} + \Vcone{1}{} + \Vcone{2}{} + \Vcone{3}{})^2 }
     =  4 \Expected{ \Vcone{0}{2} } + 12 \Expected{ \Vcone{0}{} \Vcone{1}{} } .
\end{align}
Plugging these results into Lemma \ref{lem:reduction}, we get
\begin{align}
  \Ccon{1}{3}{4}  &=  \Factor{3}{4} \cdot E_{13}^4 
             =  1536 \cdot \tfrac{3}{4} \cdot ( 2 \tfrac{1}{162} - 2 
\tfrac{1}{216} )
             =  \tfrac{32}{9}
             =  3.55\ldots ,                            \\
  \Ccon{2}{3}{4}  &=  \Factor{3}{4} \cdot E_{23}^4
             =  1536 \cdot \tfrac{1}{2} \cdot 4 \tfrac{1}{162}
             =  \tfrac{512}{27}
             =  18.96\ldots ,                           \\
  \Ccon{3}{3}{4}  &=  \Factor{3}{4} \cdot E_{33}^4
             =  1536 \cdot \tfrac{1}{8}
                     \cdot ( 4 \tfrac{1}{162} + 12 \tfrac{1}{216} )
             =  \tfrac{416}{27}
             =  15.40\ldots .
\end{align}

\ourparagraph{$4$-simplices as upper bounds.}
Here we count the critical $4$-simplices and the intervals they form
with tetrahedra, triangles, and edges as lower bounds.
For $\Ccon{3}{4}{4}$ and $\Ccon{4}{4}{4}$, we need
\begin{align}
  \tfrac{16}{5} E_{34}^4 &=  \Expected{ \Vcone{0}{} + \Vcone{1}{}
                             + \Vcone{2}{} + \Vcone{3}{} - \Vcone{4}{} }
                 =  3 \Expected{\Vcone{0}{}} ,                            \\
  16 E_{44}^4           &=  \Expected{ \Vcone{0}{} + \Vcone{1}{}
                             + \Vcone{2}{} + \Vcone{3}{} + \Vcone{4}{} }
                 =  5 \Expected{\Vcone{0}{}} .
\end{align}
We get $\Expected{\Vcone{0}{}} = \Mnt{4}{4}{1} = \tfrac{8}{81 \pi^2}$
from \eqref{eqn:moments},
and using Lemma \ref{lem:reduction}, we get
\begin{align}
  \Ccon{3}{4}{4}  &=  \Factor{4}{4} \cdot E_{34}^4 
             =  \tfrac{768 \pi^2}{5} \cdot \tfrac{5}{16} \cdot 3 \tfrac{8}{81 
\pi^2}
             =  \tfrac{128}{9}
             =  14.22\ldots ,                            \\
  \Ccon{4}{4}{4}  &=  \Factor{4}{4} \cdot E_{44}^4 
             =  \tfrac{768 \pi^2}{5} \cdot \tfrac{1}{16} \cdot 5 \tfrac{8}{81 
\pi^2}
             =  \tfrac{128}{27}
             =  4.74\ldots .
\end{align}
To avoid the complications that arise from having more than one reflection,
we compute $\Ccon{1}{4}{4}$ and $\Ccon{2}{4}{4}$
using the linear relations connecting
the Delaunay simplices with the intervals.
In particular, we get the number of tetrahedra and $4$-simplices in the
Poisson--Delaunay mosaic from the intervals as mentioned in Section \ref{sec:1};
see also \eqref{eqn:s34} and \eqref{eqn:s44}.
Since all constants other than the two sought after ones are known,
either from the above calculations or from \eqref{eqn:S44} and \eqref{eqn:S34},
this leads to a system of two linear equations:
$3 \Ccon{1}{4}{4} + 2 \Ccon{2}{4}{4} = \tfrac{737}{27}$ and
$  \Ccon{1}{4}{4} +   \Ccon{2}{4}{4} = \tfrac{346}{27}$.
Solving them, we get
\begin{align}
  \Ccon{1}{4}{4}  &=  \tfrac{  5}{ 3}  =   1.66\ldots ,                     \\
  \Ccon{2}{4}{4}  &=  \tfrac{301}{27}  =  11.14\ldots .
\end{align}
We use Lemma \ref{lem:discreteMorse} to check the number of critical simplices
and get
$\Ccon{0}{0}{4} - \Ccon{1}{1}{4} + \Ccon{2}{2}{4} - \Ccon{3}{3}{4}
                + \Ccon{4}{4}{4} = 0$, as required.

\ourparagraph{Simplices in the Poisson--Delaunay mosaic.}
Finally, we count the simplices in the Poisson--Delaunay mosaic.
Using the linear relations that connect the Delaunay simplices with the 
intervals,
we get 
\begin{align}
  \Dcon{0}{4}  &=  \Ccon{0}{0}{4}
          =  1,                                               \\
  \Dcon{1}{4}  &=  \Ccon{1}{1}{4} + \Ccon{1}{2}{4} +  \Ccon{1}{3}{4} + 
\Ccon{1}{4}{4}
          =  \tfrac{170}{9}
          =  18.88\ldots ,                                    \\
  \Dcon{2}{4}  &=  \Ccon{1}{2}{4} + \Ccon{2}{2}{4} + 2\Ccon{1}{3}{4} + 
\Ccon{2}{3}{4} + 3\Ccon{1}{4}{4} + \Ccon{2}{4}{4}
          =  \tfrac{590}{9}
          =  65.55\ldots ,                                    \\
  \Dcon{3}{4}  &=  \Ccon{1}{3}{4} + \Ccon{2}{3}{4} + \Ccon{3}{3}{4} + 
3\Ccon{1}{4}{4} + 2 \Ccon{2}{4}{4} + \Ccon{3}{4}{4}
          =  \tfrac{715}{9}
          =  79.44\ldots ,               \label{eqn:s34}      \\
  \Dcon{4}{4}  &=  \Ccon{1}{4}{4} + \Ccon{2}{4}{4} + \Ccon{3}{4}{4} + 
\Ccon{4}{4}{4}
          =  \tfrac{286}{9}
          =  31.77\ldots .               \label{eqn:s44}
\end{align}
This completes the justification of the numbers in Tables \ref{tbl:constants}
and \ref{tbl:PD-constants}.
We note that we did not use the Euler Relations to derive any of the constants.
We can therefore use it to check whether the computations are possibly correct.
Indeed, we get
$\Dcon{0}{4} - \Dcon{1}{4} + \Dcon{2}{4} - \Dcon{3}{4} + \Dcon{4}{4} = 0$,
as required.

\section{Discussion}
\label{sec:6}

Using a Poisson point process to sample a random set of points in $\Rspace^n$,
we study the radius function on the Poisson--Delaunay mosaic \cite{BaEd15}.
Our main result are integral expressions for the expected numbers
of critical simplices and intervals of this generalized discrete Morse function
that depend on the maximum allowed radius.
This work suggests a number of open questions.

\begin{enumerate}\denselist
  \item[1.]  We have concrete expressions for the constants that show up
    in the expectations in dimensions $n = 2, 3, 4$.
    With the exception of dimensions $0, n-1, n$ \cite{ScWe08},
    even the expected numbers of simplices in the Poisson--Delaunay mosaic
    are currently not known beyond four dimensions.
    What is the asymptotic behavior of the constants
    $\Ccon{\ell}{k}{n}$ and $\Dcon{k}{n}$ as $n$ goes to infinity?
  \item[2.] Can the results be extended to weighted Delaunay mosaic
    as defined in \cite{Ede01}?
    We suggest that a natural model of such a mosaic
    is the nerve of an affine slice of a Poisson-Voronoi diagram in $\Rspace^n$.
    Observe however that the implied distribution of the radii depends
    on the co-dimension of the slice.
  \item[3.] Can the results be extended to the Betti numbers
    and the framework of persistent homology; see e.g.\ \cite{BKS15}?
    Indeed, the intervals of size larger than $1$ correspond to
    $0$-persistent pairs,
    and it is natural to ask similar questions about the pairs with
    positive persistence.
\end{enumerate}

\noindent
We finally mention that it is interesting to ask the same questions for
the \v{C}ech complexes of a Poisson point process in $\Rspace^n$.
Mapping each simplex to the radius of the smallest enclosing sphere,
we get again a generalized discrete Morse function; see \cite{BaEd15}.
The critical simplices are the same as for the radius function of the
Poisson--Delaunay mosaic, but there is a much richer structure of intervals,
which can be analyzed with the methods of this paper.


\newpage
\appendix
\section{Boundary Effect}
\label{app:A}

Recall that $K_0$ is the nerve of the Voronoi diagram restricted to $\Region$,
and $K_1 \subseteq K_0$ contains all Delaunay simplices whose smallest
empty circumspheres have the center inside $\Region$.
In this appendix, we show that the difference between $K_0$ and $K_1$ is small
when $\Region$ is a ball.

\ourparagraph{Big spheres.}
We need an auxiliary lemma implying that only a vanishing fraction
of the $n$-simplices in the Poisson--Delaunay mosaic have circumspheres with
radii larger than some positive threshold.
To simplify the discussion, we call the closed ball bounded by the
circumsphere of an $n$-simplex its \emph{circumball}.
Letting $\Hspace \subseteq \Rspace^n$ be bounded and $r_0 > 0$,
we write $\# (\Hspace, r_0)$ for the number of $n$-simplices in the
Poisson--Delaunay mosaic whose circumspheres have the center in $\Hspace$ and
the radius exceeding $r_0$.
\begin{lemma}[Big Spheres]
  \label{lem:bigspheres}
  There exist positive constants $c, \alpha, \beta$, all depending only on $n$,
  such that for any bounded Borel region $\Hspace \subseteq \Rspace^n$ and any fixed 
$r_0 > 0$,
  $\Expected{\#(\Hspace, r_0)} \leq c \norm{\Hspace} e^{-\alpha r_0^\beta}$.
\end{lemma}
\ourproof
  Arguing as in the proof of Lemma \ref{lem:counting} with the only difference
  that $z$ is now integrated over $\Hspace$ and $r$ from $r_0$ to infinity, we 
can write
  \begin{align}
    \Expected{\#(\Hspace, r_0)}  &=  c_0 \norm{\Hspace} \int_{r=r_0}^{\infty}
                               e^{- \density r^n \nu_n} r^{n^2-1} \diff r 
				  = c_0' \norm{\Hspace} \Gamma(n, \density 
r_0^n 
\nu_n),
      \label{eqn:bs2}
  \end{align}
  in which $c_0$ and $c_0'$ are constants that depends only on $n$,
  and $\Gamma(k, x) = \Gamma(k) - \gamma(k,x)$ is the
  upper incomplete Gamma function.
  Noticing that
  $\Gamma(n, \density r_0^n \nu_n) = \littleoh{e^{-0.9 \density r_0^n \nu_n}}$,
  see for example \cite{Olv97}, completes the proof.
\eop

\ourparagraph{Size of boundary.}
We are now ready to give an upper bound on the number of
simplices in $K_0$ that are not in $K_1$, which we need for
bound \eqref{eqn:Euler-D} on the Euler characteristic of $K_1$.
Every simplex $Q \in K_0 \setminus K_1$ corresponds to an intersection
of Voronoi domains, $\Voronoi{Q}$, that has points inside as well as outside 
$\Region$.
Let $x \in \Voronoi{Q} \cap \Region$ and $y \in \Voronoi{Q} \setminus \Region$.
We argue that both points are contained in 
the union of circumballs of the $n$-simplices that share $Q$.
Indeed, all these circumballs contain all points of $Q$,
and for each $q \in Q$ there is a vertex of $\Voronoi{Q}$ that is closer
to $x$ than to $q$, so the circumball centered at this vertex contains $x$.
The same argument applies to $y$.
Since the union contains points on both sides of $\partial \Region$,
at least one of these circumballs has a non-empty
intersection with $\partial \Region$.

Writing $\# (\partial \Region)$ for the number of $n$-simplices whose 
circumballs
have a non-empty intersection with $\partial \Region$,
we prove that it grows slower than the number of $n$-simplices
whose circumballs are centered inside $\Region$. The discussion above implies 
that
$|K_0 \setminus K_1| < 2^{n+1}\# (\partial \Region)$,
so to have $|K_0 \setminus K_1| = \littleoh{\density \norm{\Region}}$,
it is enough to prove the following.
\begin{lemma}[Boundary Size]
  \label{lem:boundarysize}
  Let $X$ be a Poisson point process with density $\density$
  in $\Rspace^n$. Let $\Omega = \Bspace(R)$ be a ball of radius $R$ centered at the origin.
  Then $\Expected{\# (\partial \Region)} = \littleoh{1} \density 
\norm{\Region}$.
\end{lemma}
\ourproof
  Without loss of generality assume $\density = 1$.
  It suffices to count the $n$-simplices with circumcenters outside $\Region$
  and to prove that the number of such $n$-simplices
  whose circumballs intersect $\partial \Region$ is $\bigoh{R^{n-1+\dd}}$.
  Assume $R > 1$, fix $0 < \dd < 1$ and
  let $\Aspace$ be the set of points at distance at most $R^\dd$ from $\partial 
\Region$.
  For a ball with center $z$ outside $\Region$ to intersect $\Aspace$,
  one of the following must happen:
  \begin{enumerate}\denselist
    \item[1.]  $z \in \Aspace$;
    \item[2.]  $z \in \Bspace (2R) \setminus \Aspace$ and its radius exceeds 
$R^\dd$;
    \item[3.]  $z \not\in \Bspace (2R)$ and its radius exceeds $R$.
  \end{enumerate}
  As proved in \cite{ScWe08} and reproved in this paper,
  the expected number of $n$-simplices in $\Delaunay{S}$ with center
  in $\Aspace$ is $\bigoh{\norm{\Aspace}} = \bigoh{R^{n-1+\dd}}$.
  This settles Case 1.
  Applying Lemma \ref{lem:bigspheres}, we see that the expected number of
  $n$-simplices with center in $\Bspace (2R)$
  and radius larger than $R^\dd$ is
  $\bigoh{R^n e^{- \alpha R^{\dd\beta} }}$,
  in which $\alpha$ and $\beta$ are positive constants.
  This settles Case 2.
  Finally, we decompose the complement of $\Bspace (2R)$
  into annuli of the form
  $\Hspace_i = \Bspace (iR+2R) \setminus \Bspace (iR+R)$,
  for $i \geq 1$.
  To intersect $\Region = \Bspace (R)$, a ball centered inside $\Hspace_i$
  must have radius exceeding $i R$.
  Writing $\Hspace = \bigcup_{i=1}^\infty \Hspace_i$ for the union of annuli
  and $\# (\Hspace, \Region)$ for the number of $n$-simplices
  with circumcenter in $\Hspace$ whose circumball intersects $\Region$,
  we get an upper bound on the expected number:
  \begin{align}
    \Expected{\#(\Hspace, \Region)}
       &\leq  \sum_{i=1}^{\infty} \Expected{\#(\Hspace_i, i r)}    
        \label{eqn:outer-expectation1}                                    \\
       &\leq  \sum_{i=1}^{\infty} c \norm{\Hspace} e^{-\alpha (ir)^\beta}
        \label{eqn:outer-expectation2}                                    \\
      &\leq  c' R^n e^{-\alpha R^\beta}
             \sum_{i=0}^{\infty} i^n e^{-\alpha i^\beta},
        \label{eqn:outer-expectation3}
  \end{align}
  where we use Lemma \ref{lem:bigspheres} to get \eqref{eqn:outer-expectation2},
  and $\norm{\Hspace_i} = \bigoh{i^n R^n}$
  as well as $\alpha (iR)^\beta \geq \alpha i^\beta$
  to get \eqref{eqn:outer-expectation3}.
  Since the last sum converges, we get
  $\Expected{\#(\Hspace, \Region)} = \bigoh{R^n e^{-\alpha R^\beta}}$,
  which settles Case 3.
\ignore
{
  Taking $\dd = \min \{ \tfrac{\ee}{2}, 1\}$ and $\gamma > 0$,
  summing up results from all three cases, and applying Chebyshev's inequality 
we get
  \begin{align}
    \Probable{}{\tfrac{\# (\partial \Region)}{\norm{\partial\Region}^{1+\ee}} > 
\gamma}
      &< \tfrac{1}{\gamma} \tfrac{\Expected{\# (\partial \Region)}}
                                                   
{\norm{\partial\Region}^{1+\ee}} \\
      &= \tfrac{1}{\gamma}
                           \littleoh{\tfrac{1}{\norm{\partial\Region}^{\ee - 
\dd}}},
      \label{eqn:cheb}
  \end{align}
  which goes to $0$ and implies
  $\# (\partial \Region) = \littleoh{\norm{\partial\Region}^{1+\ee}}$,
  as claimed.
}
\eop

\noindent {\sc Remarks.}
Besides $| K_0 \setminus K_1 | = \littleoh{1} \density \norm{\Region}$,
Lemma \ref{lem:boundarysize} implies that the number of vertices of $K_0$
outside $\Region$ is $\littleoh{1} \density \norm{\Region}$.

Actually, we have proved that for any $\ee > 0$,
$\Expected{\# (\partial \Region)} = \littleoh{1} \density 
\norm{\partial\Region}^{1+\ee}$.
Also, one can apply the Markov's inequality to show that the convergence
happens in probability.


\begin{thebibliography}{15}

\bibitem{BaEd15}
{\sc U.\ Bauer and H.\ Edelsbrunner.}
The Morse theory of \v{C}ech and Delaunay complexes.
\emph{Trans.\ Amer.\ Math.\ Soc.}, to appear.

\bibitem{BaLa09}
{\sc V. Baumstark and G. Last.}
Gamma distributions for stationary Poisson flat processes.
\emph{Adv.\ Appl.\ Prob.} {\bf 41} (2009), 911--939.

\bibitem{BoAd14}
{\sc O.\ Bobrowski and R.\ Adler.}
Distance functions, critical points, and topology for some random complexes.
\emph{Homology, Homotopy and Applications} {\bf 16} (2014), 311--344.

\bibitem{BKS15}
{\sc O.\ Bobrowski, M.\ Kahle and P.\ Skraba.}
Maximally persistent cycles in random geometric complexes.
arXiv:1509.04347, 2015.

\bibitem{BoWe15}
{\sc O.\ Bobrowski and S.\ Weinberger.}
On the vanishing of homology in random \v{C}ech complexes.
arXiv:1507.06945v1, 2015.


\bibitem{Car09}
{\sc G.\ Carlsson.}
Topology and data.
\emph{Bull.\ Amer.\ Math.} {\bf 46} (2009), 255--308.

\bibitem{Che14}
{\sc N.\ Chenavier.}
A general study of extremes of stationary tessellations with applications.
\emph{Stoch.\ Process.\ Appl.} {\bf 124} (2014), 2917--2953.

\bibitem{DFRV14}
{\sc L.\ Decreusefond, E.\ Ferraz, H.\ Randriam and A.\ Vergne.}
Simplicial homology of random configurations.
\emph{Adv.\ Appl.\ Prob.} {\bf 46} (2014), 1--23.

\bibitem{Dey11}
{\sc T.K.\ Dey.}
\emph{Curve and Surface Reconstruction.  Algorithms with Mathematical Analysis.}
Cambridge Univ.\ Press, Cambridge, England, 2011.

\bibitem{Ede01}
{\sc H.\ Edelsbrunner.}
\emph{Geometry and Topology for Mesh Generation.}
Cambridge Univ.\ Press, Cambridge, England, 2001.

\bibitem{Ede03}
{\sc H.\ Edelsbrunner.}
Surface reconstruction by wrapping finite sets of points in space.
\emph{Discrete and Computational Geometry.  The Goodman-Pollack Festschrift},
  379--404, eds.\ B.\ Aronov, S.\ Basu, J.\ Pack and M.\ Sharir,
  Springer-Verlag, 2003.

\bibitem{EdHa10}
{\sc H.\ Edelsbrunner and J.L.\ Harer.}
\emph{Computational Topology. An Introduction.}
Amer.\ Math.\ Soc., Providence, Rhode Island, 2010.

\bibitem{EdMu94}
{\sc H.\ Edelsbrunner and E.P.\ M\"{u}cke.}
Three-dimensional alpha shapes.
\emph{ACM Trans.\ Graphics} {\bf 13} (1994), 43--72.

\bibitem{For98}
{\sc R.\ Forman.}
Morse theory for cell complexes.
\emph{Adv.\ Math.} {\bf 134} (1998), 90--145.

\bibitem{Fre09}
{\sc R.\ Freij.}
Equivariant discrete Morse theory.
\emph{Discrete Math.} {\bf 309} (2009), 3821--3829.

\bibitem{Kah11}
{\sc M.\ Kahle.}
Random geometric complexes.
\emph{Discrete Comput.\ Geom.} {\bf 45} (2011), 553--573.

\bibitem{Kah14}
{\sc M.\ Kahle.}
Topology of random simplicial complexes:  a survey.
\emph{AMS Contemp. Math}. {\bf 620} (2014), 201--222.

\bibitem{Kin93}
{\sc J.F.C.\ Kingman.}
\emph{Poisson Processes.}
Oxford Univ.\ Press, Oxford, England, 1993.

\bibitem{Mil69}
{\sc R.E.\ Miles.}
Poisson flats in Euclidean spaces.
  Part I:  a finite number of random uniform flats.
\emph{Adv.\ Appl.\ Prob.} {\bf 1} (1969), 211--237.

\bibitem{Mil70}
{\sc R.E.\ Miles.}
On the homogeneous planar Poisson point process.
\emph{Math.\ Biosci.} {\bf 6} (1970), 85--127.

\bibitem{Mil71}
{\sc R.E.\ Miles.}
Isotropic random simplices.
\emph{Adv.\ Appl.\ Prob.} {\bf 3} (1971), 353--382.

\bibitem{Mil74}
{\sc R.E.\ Miles.}
A synopsis of `Poisson flats in Euclidean spaces'. 
In \emph{Stochastic Geometry}, eds.: E.F.\ Harding and D.G.\ Kendall,
  John Wiley, New York, 202--227, 1974.

\bibitem{Mol89}
{\sc J.\ M{\o}ller.}
Random tessellations in $\Rspace^d$.
\emph{Adv.\ Appl.\ Prob.} {\bf 21} (1989), 37--73.

\bibitem{Olv97}
{\sc F.W.J.\ Olver.}
\emph{Asymptotics and Special Functions.}
A.K.\ Peters, Wellesley, Massachusetts, 1997.

\bibitem{ScWe08}
{\sc R.\ Schneider and W.\ Weil.}
\emph{Stochastic and Integral Geometry.}
Springer, Berlin, Germany, 2008.

\bibitem{Wen62}
{\sc J.G.\ Wendel.}
A problem in geometric probability.
\emph{Math.\ Scand.} {\bf 11} (1962), 109--111.

\end{thebibliography}
\end{document}